\documentclass{amsart}

\usepackage{latexsym}

\usepackage{amsfonts,amssymb,euscript,epsfig}

\usepackage[all]{xy}

\newtheorem{theo}{Th\'eor\`eme}[section]
\newtheorem{defi}[theo]{D\'efinition}

\newtheorem{lem}[theo]{Lemme}
\newtheorem{prop}[theo]{Proposition}

\newtheorem{rem}[theo]{Remarque}

\newcommand{\ggot}{\ensuremath{\mathfrak{g}}}

\newcommand{\kgot}{\ensuremath{\mathfrak{k}}}

\newcommand{\tgot}{\ensuremath{\mathfrak{t}}}
\newcommand{\Rgot}{\ensuremath{\mathfrak{R}}}

\newcommand{\Acal}{\ensuremath{\mathcal{A}}}
\newcommand{\Bcal}{\ensuremath{\mathcal{B}}}
\newcommand{\Ccal}{\ensuremath{\mathcal{C}}}
\newcommand{\Dcal}{\ensuremath{\mathcal{D}}}
\newcommand{\Ecal}{\ensuremath{\mathcal{E}}}

\newcommand{\Fcal}{\ensuremath{\mathcal{F}}}
\newcommand{\Hcal}{\ensuremath{\mathcal{H}}}

\newcommand{\Lcal}{\ensuremath{\mathcal{L}}}
\newcommand{\Qcal}{\ensuremath{\mathcal{Q}}}
\newcommand{\Mcal}{\ensuremath{\mathcal{M}}}
\newcommand{\Ncal}{\ensuremath{\mathcal{N}}}
\newcommand{\Ocal}{\ensuremath{\mathcal{O}}}

\newcommand{\Rcal}{\ensuremath{\mathcal{R}}}

\newcommand{\Ycal}{\ensuremath{\mathcal{Y}}}
\newcommand{\Zcal}{\ensuremath{\mathcal{Z}}}
\newcommand{\Ucal}{\ensuremath{\mathcal{U}}}

\newcommand{\Z}{\ensuremath{\mathbb{Z}}}
\newcommand{\C}{\ensuremath{\mathbb{C}}}
\newcommand{\Lfibre}{\ensuremath{\mathbb{L}}}

\newcommand{\R}{\ensuremath{\mathbb{R}}}
\newcommand{\N}{\ensuremath{\mathbb{N}}}
\newcommand{\tore}{\ensuremath{\mathbb{T}}}

\newcommand{\esp}{\ensuremath{\varepsilon}}
\newcommand{\f}{\ensuremath{\mathcal{C}^{\infty}}}
\newcommand{\fgene}{\ensuremath{\mathcal{C}^{-\infty}}}

\newcommand{\croc}{\ensuremath{\hookrightarrow}}

\newcommand{\indice}{\ensuremath{\hbox{\rm Indice}}}
\newcommand{\ch}{\ensuremath{\hbox{\rm H}}}
\newcommand{\lde}{\ensuremath{\hbox{\rm L}^2}}
\newcommand{\loc}{\ensuremath{\hbox{\rm P}}}
\newcommand{\p}{\ensuremath{\Delta}}
\newcommand{\mm}{\ensuremath{\hbox{\rm m}}}
\newcommand{\mdh}{\ensuremath{\hbox{\rm DH}}}
\newcommand{\noyau}{\ensuremath{\hbox{\rm Ker}}}
\newcommand{\kir}{\ensuremath{\hbox{\rm Kir}}}
\newcommand{\kv}{\ensuremath{\hbox{\rm kv}}}
\newcommand{\crit}{\ensuremath{\hbox{\rm Cr}}}
\newcommand{\Vol}{\ensuremath{\hbox{\rm vol}}}
\newcommand{\T}{\ensuremath{\hbox{\bf T}}}

\newcommand{\K}{\ensuremath{\hbox{\bf K}}}
\newcommand{\Eul}{\ensuremath{\hbox{\rm Eul}}}
\newcommand{\Char}{\ensuremath{\hbox{\rm Char}}}

\newcommand{\Thom}{\ensuremath{\hbox{\rm Thom}}}

\def \spin {{\small\rm spin}} 
 
\def \hol {{\small\rm hol}}
\def \spinc {{\rm Spin}^{c}} 
\def \so {{\rm so}} 
 
\def \det {{\rm det}}

\def \indB {{\rm Ind}^{K}_{K_{\beta}}} 
\def \indt {{\rm Ind}^{^K}_{_T}} 
\def \indh {{\rm Ind}^{^K}_{_H}}

\def \HolB {{\rm Hol}^{K}_{K_{\beta}}}

\begin{document}


\title{Cohomologie \'equivariante et quantification g\'eom\'etrique}

\vspace{1cm}

\author{Paul-\'Emile PARADAN}

\maketitle

{\center UMR 5582, Institut Fourier, B.P. 74, 38402,
Saint-Martin-d'H\`eres
cedex, France\\
e-mail: Paul-Emile.Paradan@ujf-grenoble.fr\\
}

\begin{abstract}
Ce m\'emoire est le texte de mon Habilitation \`a Diriger des Recherches. 
Dans celui-ci je rappelle les techniques que j'ai mises en oeuvre pour r\'ealiser 
le programme {\em de localisation non-ab\'elienne} de Witten, et les 
r\'esultats qui en d\'ecoulent. Ces travaux de recherches concernent les diff\'erentes 
th\'eories cohomologiques associ\'ees aux actions de groupes de Lie compacts sur des 
vari\'et\'es diff\'erentiables: cohomologie \'equivariante, K-th\'eorie \'equivariante, 
et la th\'eorie des op\'erateurs transversalement elliptiques. 
\end{abstract}

{\small
\tableofcontents
}

\section{Cohomologie \'equivariante et localisation}

La ``cohomologie \'equivariante'' est n\'ee dans les ann\'ees 50 apr\`es
les travaux de  Borel et de H. Cartan \cite{Cartan}. Une version topologique du
th\'eor\`eme de localisation appara\^\i t dans les travaux de Borel
\cite{Borel} et Quillen \cite{Quillen}, et en K-th\'eorie
\'equivariante dans ceux de Atiyah et Segal \cite{Segal68}. Il faudra
attendre une quinzaine d'ann\'ees pour que, sur l'impulsion de la
formule de Duistermaat-Heckman, on obtienne le th\'eor\`eme de
localisation de Berline-Vergne et d'Atiyah-Bott \cite{Atiyah-Bott,B-Vergne1}.  
Berline et Vergne obtiennent une localisation de {\em l'int\'egrale} d'une forme
\'equivariante. Un peu plus tard, Atiyah et Bott raffinent cette
localisation au niveau de la cohomologie. Ces travaux s'obtiennent
dans le cas d'un groupe {\em ab\'elien}. En 1992, Witten propose une
localisation {\em non-ab\'elienne} dans le cadre Hamiltonien \cite{Witten}. 
Mon travail dans ce domaine a \'et\'e en grande partie consacr\'ee \`a la
r\'ealisation du programme de Witten.

Dans une premi\`ere partie je rappelle succinctement le mod\`ele de
Cartan de la cohomologie \'equivariante.  Avant de d\'ecrire les
diff\'erentes localisations dont j'ai parl\'e plus haut, j'\'evoque le
travail pr\'ecurseur de Bott \cite{Bott67} sur les nombres
caract\'eristiques: on verra que la m\'ethode de Bott contient en
substance le proc\'ed\'e de localisation de Berline-Vergne.

\subsection{Mod\`ele de Cartan}\label{subsec:modele-cartan}
Consid\'erons un groupe de Lie compact connexe $K$, d'alg\`ebre de Lie
$\kgot$, agissant de mani\`ere $\f$ sur une vari\'et\'e
diff\'erentielle $M$. On note $\Acal(M)$ l'alg\`ebre sur $\C$ des
formes diff\'erentielles $\f$ et $d$ la d\'erivation de de Rham. Si $\xi$
est un champ de vecteurs sur $M$, on note $c(\xi):\Acal(M)\to\Acal(M)$
la contraction par $\xi$. L'action de $K$ sur $M$
d\'etermine un morphisme $X\to X_{M}$ de $\kgot$ dans l'alg\`ebre des
champs de vecteurs de $M$.

On consid\`ere l'espace des applications $K$-\'equivariantes
$\kgot\to\Acal(M),\ X\mapsto\eta(X)$, muni de la d\'erivation $D$
\begin{equation}
  \label{eq:def-D}
  (D\eta)(X):= (d-c(X_{M}))(\eta(X)),\ X\in\kgot.
\end{equation}
Comme $D^2=0$, on peut consid\'erer l'espace de cohomologie 
$\noyau D/{\rm Im} D$.

Le mod\`ele de Cartan \cite{B-G-V,Cartan,Duflo-Vergne2,Guill-Stern00} consid\`ere des 
applications {\em  polynomiales}, et l'espace de cohomologie associ\'e est not\'e
$\Hcal_{K}^{*}(M)$. Cet espace est muni naturellement d'une structure
d'alg\`ebre sur $\Hcal_{K}^{*}(\cdot)=S(\kgot^*)^K$, l'alg\`ebre 
(sur $\C$) des applications polynomiales $K$-invariantes sur $\kgot$.

On peut aussi consid\'erer des applications $X\mapsto\eta(X)$ qui
sont $\f$ , et on obtient comme espace de cohomologie l'alg\`ebre
$\Hcal_{K}^{\infty}(M)$. S. Kumar et Vergne \cite{Kumar-Vergne} ont
\'etudi\'e l'espace de cohomologie $\Hcal_{K}^{-\infty}(M)$ obtenu
en consid\'erant des applications $X\mapsto\eta(X)$ qui sont $\fgene$.
Rappelons sa construction et les diff\'erentes notations.

Soit $\fgene(\kgot,\Acal(M))$ l'espace des fonctions g\'en\'eralis\'ees 
sur $\kgot$ \`a valeurs dans $\Acal(M)$. C'est, par d\'efinition, l'espace 
des applications $\C$-lin\'eaires continues de l'espace des densit\'es $\f$ \`a
support compact de $\kgot$ dans $\Acal(M)$. L'image de la densit\'e 
$\phi(X)dX$ par $\eta\in\fgene(\kgot,\Acal(M))$ est une
forme diff\'erentielle sur $M$ not\'ee $<\eta(X),\phi(X)dX>_\kgot$.
La diff\'erentielle $D$ d\'efinie sur $\f(\kgot,\Acal(M))$ par 
(\ref{eq:def-D}) se prolonge \`a $\fgene(\kgot,\Acal(M))$, et on montre 
que $D^2=0$ sur le sous-espace $\fgene(\kgot,\Acal(M))^K$ des \'el\'ements 
$K$-invariants \cite{Kumar-Vergne}. L'espace de cohomologie associ\'e est appel\'e
la cohomologie $K$-\'equivariante \`a coefficients
g\'en\'eralis\'es de $M$, et est not\'e $\Hcal_{K}^{-\infty}(M)$.
 On remarque que ce dernier espace ne poss\`ede plus de structure
multiplicative, mais est muni n\'eanmoins d'une structure de
$\Hcal_{K}^{*}(M)$-module.  

\medskip

Les espaces de cohomologie $\Hcal_{K}^{-}(M)$, pour 
$-\in\{*,\infty,-\infty\}$, ont des bonnes propri\'et\'es fonctorielles. Si
$i:N\croc M$ est une sous-vari\'et\'e $K$-stable, on a un morphisme de
{\em restriction} $i^*:\Hcal_{K}^{-}(M) \to\Hcal_{K}^{-}(N)$. Si de
plus le fibr\'e normal $\Ncal$ de $N$ par rapport \`a $M$ est {\em
  orient\'e}, nous avons un morphisme `{\em image directe}'
$i_*:\Hcal_{K}^{-}(N) \to\Hcal_{K}^{-}(M)$. Leur compos\'ee $i^*\circ
i_*$ est le morphisme de multiplication par la classe d'Euler
equivariante, $\Eul(\Ncal)\in\Hcal_K^*(N)$. Historiquement, la
construction du morphisme image directe $i_*$ dans le mod\`ele de
Cartan remonte \`a l'article de Mathai-Quillen \cite{Mathai-Quillen},
o\`u ils explicitent un repr\'esentant de la classe de Thom
\'equivariante.

D'autre part, si $P\to M$ est une fibration $K$-\'equivariante {\em
  orient\'ee}, on a un morphisme d'{\em int\'egration} le long des
fibres $\int_{P/M}:\Hcal_{K}^{-}(P)\to\Hcal_{K}^{-}(M)$. En
particulier si $M$ est orient\'ee, on a un morphisme\footnote{
$\mathcal{C}^{*}(\kgot)^K=S(\kgot^*)^K$.} d'int\'egration
$\int_M:\Hcal_{K}^{-}(P)\to\mathcal{C}^{-}(\kgot)^K$.  

\medskip

Si on consid\`ere un sous-groupe de Lie $H\croc K$, on a au niveau des
coefficients un morphisme de {\em restriction}
$\Hcal_{K}^{-}(M)\to\Hcal_{H}^{-}(M),\ \alpha\mapsto \alpha\vert_H$
pour $-\in\{*,\infty\}$.  Pour les coefficients g\'en\'eralis\'es, les
choses se passent de mani\`ere {\em duale}.  Kumar et Vergne
d\'efinissent, lorsque $H$ et $K$ ont le m\^eme rang, un morphisme
d'{\em induction}
$\indh:\Hcal_{H}^{-\infty}(M)\to\Hcal_{K}^{-\infty}(M)$ qui poss\`ede
la bonne compatibilit\'e par rapport au morphisme de restriction et
les structures de module: pour $\alpha\in\Hcal_{K}^{*}(M)$ et
$\beta\in\Hcal_{H}^{-\infty}(M)$ on a
$\alpha\cdot\indh(\beta)=\indh(\alpha\vert_H\cdot\beta)$.

\bigskip

{\em Un exemple fondamental de forme \'equivariante \`a coefficients
  g\'en\'eralis\'es}.  Supposons que le groupe de Lie compact $K$
agisse librement sur $M$. Notons $\pi:M\to B$ le fibr\'e principal
correspondant. Soit $\sigma\in(\Acal^1(M)\otimes\kgot)^K$ une
$1$-forme de connexion, et $\Omega=d\sigma+\frac{1}{2}[\sigma,\sigma]$ 
sa courbure. Nous avons l'isomorphisme de Chern-Weil
$$
{\rm cw}:\Hcal_K^*(M)\to\Hcal^*(B) \ .
$$

Consid\'erons sur $M$ la forme \'equivariante \`a coefficients
g\'en\'eralis\'es $\delta(X-\Omega)$, qui est d\'efinie par la
relation: $<\delta(X-\Omega),\phi(X)dX>_\kgot= \phi(\Omega)\Vol(K,dX)$, o\`u
$\Vol(K,dX)$ est le volume de $K$ pour la mesure de Haar compatible 
avec $dX$. Le terme $\phi(\Omega)$ est la valeur
de l'op\'erateur diff\'erentiel $e^{\Omega(\frac{\partial}{\partial X}|_0)}$ contre
la fonction $\phi$. Soient $\sigma_k$ les composantes de la $1$-forme de connexion 
relativement \`a une base $E_1,\ldots,E_r$ de $\kgot$. S. Kumar et Vergne ont introduit
la forme \'equivariante ferm\'ee
$\delta(X-\Omega)\frac{\sigma_1\wedge\ldots\wedge\sigma_r}{\Vol(K)}$, avec laquelle
ils d\'eterminent l'isomorphisme
\begin{eqnarray}\label{eq-kv}
{\rm kv}:\Hcal^*(B)&\to&\Hcal_K^{-\infty}(M)\\  \nonumber
\eta&\mapsto&\pi^*(\eta)\delta(X-\Omega)\frac{\sigma_1\wedge\ldots\wedge\sigma_r}
{\Vol(K)}    \ .
\end{eqnarray}
Ici $\Vol(K)$ est le volume de $K$ pour la mesure de Haar compatible avec
la base $E_1,\ldots,E_r$. Les isomorphismes `cw' et `kv' sont de plus compatibles:
$\eta\wedge{\rm kv}(\gamma)={\rm kv}({\rm cw}(\eta)\wedge\gamma)$ pour
$\eta\in\Hcal_K^{*}(M)$ et $\gamma\in\Hcal^*(B)$. On remarque dans ce cas 
que $\kv\circ{\rm cw}$ est un isomorphisme de $\Hcal_K^*(M)$ sur
$\Hcal_K^{-\infty}(M)$, tandis que le morphisme naturel $\Hcal_K^*(M)
\to\Hcal_K^{-\infty}(M)$ d'extension des coefficients est le {\em
  morphisme nul}.

{\bf Convention}: Pour une fonction g\'en\'eralis\'ee $f\in\fgene(\kgot)$ support\'ee
en $0$, on parle de {\em sa multiplicit\'e par rapport \`a la masse de Dirac en $0$}.
C'est la quantit\'e $<f(X),\phi(X)dX>_\kgot\in\C$, o\`u $\phi$ est une fonction 
$\f$ \'egale \`a $1$ au voisinage de $0\in\kgot$, et $dX$ est normalis\'ee par la 
condition $\Vol(K,dX)=1$.

Ici, pour tout $\eta\in\Hcal^*(B)$, la fonction g\'en\'eralis\'ee $\int_M\kv(\eta)$ 
est support\'ee en $0$, et sa multiplicit\'e par rapport \`a la masse de Dirac en $0$ est
\'egale \`a $\int_B\eta$.

\subsection{Champs de vecteurs et nombres caract\'eristiques}

Je reprends ici le titre original de l'article de Bott \cite{Bott67}. 
Un th\'eor\`eme de Hopf affirme que la
caract\'eristique d'Euler d'une vari\'et\'e compacte est \'egal au
nombre de z\'eros d'un champ de vecteurs (compt\'es judicieusement).
Dans cet article Bott d\'emontre qu'un principe similaire s'applique
aux {\em nombres de Pontrijagin}.

Soit $M$ une vari\'et\'e riemanienne compacte orient\'ee de dimension
$2n$. Les nombres de Pontrijagin s'expriment au moyen des classes de
Chern de l'espace tangent complexifi\'e $\T M\otimes\C$. Il sont aussi
d\'etermin\'es par le proc\'ed\'e suivant. Soit $\nabla$ la connexion
de Levi-Civita sur $M$, et $R=\nabla^2$ le tenseur de courbure
associ\'e. \`A chaque application polynomiale homog\`ene\footnote{Ici $\so(2n)$ 
est l'alg\`ebre de Lie de $SO(2n)$.} $\phi:\so(2n)\to \C$
invariante par rapport \`a l'action adjointe de $SO(2n)$, on peut
associer la forme diff\'erentielle ferm\'ee $\phi(R)$ ainsi que son
int\'egrale $\phi[M]:=\int_M\phi(R)$. Si le degr\'e de l'application 
polynomiale homog\`ene $\phi$ est \'egal \`a  $n$, $\phi[M]$ est un nombre 
caract\'eristique; sinon $\phi[M]=0$.

Consid\'erons un champ de vecteurs de {\em Killing} $V$ sur $M$, tel
que l'ensemble des z\'eros de $V$, not\'e $M^V$, est discret. Dans cet
article, Bott montre que les int\'egrales $\phi[M]$ se localisent sur
$M^V$. Le champ de vecteurs $V$ d\'etermine un action $\Lcal(V)$ sur les sections 
du fibr\'e tangent $\T M$. Pour chaque $p\in M^V$, cette action se sp\'ecialise
en un endomorphisme $\Lcal(V)_p$ de l'espace tangent $\T_p(M)$: l'endomorphisme 
est de plus antisym\'etrique car $V$ est de Killing.  D'une part, l'orientation 
de $M$ permet de d\'efinir une racine carr\'ee de $\det(\Lcal(V)_p)$: le Pfaffien
$\det^{1/2}(\Lcal(V)_p)$. D'autre part le polyn\^ome invariant $\phi$
d\'etermine $\phi(\Lcal(V)_p)\in \C$. La formule de Bott s'\'enonce
ainsi. Pour tout polyn\^ome $\phi$ de degr\'e {\em inf\'erieur o\`u \'egal} \`a
$n$, on a
\begin{equation}
  \label{eq:th-bott}
(-2\pi)^n \sum_{p\in M^V}\frac{\phi(\Lcal(V)_p)}{\det^{1/2}(\Lcal(V)_p)}=\phi[M]\ . 
\end{equation}

\medskip

Voici une des \'etapes clef de la d\'emonstration de Bott. Pour chaque
$\phi$ de degr\'e {\em inf\'erieur o\`u \'egal} \`a $n$, il explicite une forme
diff\'erentielle de degr\'e $2n-1$ sur l'ouvert $M-M^V$, not\'ee
$\eta_{\phi}$, telle que
\begin{equation}\label{eq-bott-1}
[\phi(R)]^{\rm max}=d\eta_{\phi} \quad {\rm sur}\quad M-M^V\ . 
\end{equation}
o\`u $[-]^{\rm max}$ d\'esigne la composante de degr\'e maximal $2n$.  Soit
$M_{\esp}$ le compl\'ementaire d'un $\esp$-voisinage de $M^V$. Alors
$\phi[M]=\lim_{\esp\to 0}\int_{M_\esp}[\phi(R)]^{\rm max}$. Mais sur
$M_\esp$, le th\'eor\`eme de Stokes combin\'e avec (\ref{eq-bott-1})
donne $\int_{M_\esp}[\phi(R)]^{\rm max} =
\int_{M_\esp}d\eta_\phi=\int_{\partial M_\esp}\eta_\phi$. On obtient
alors $\phi[M]=\lim_{\esp\to 0}\int_{\partial M_\esp}\eta_\phi$, et
l'expression (\ref{eq:th-bott}) se d\'emontre apr\`es le calcul
(local) de $\int_{\partial M_\esp}\eta_\phi$.

Nous finissons cette section avec la construction de cette forme
$\eta_\phi$. On verra l\`a les pr\'emices, quinze ans auparavant, de
la m\'ethode de localisation en cohomologie \'equivariante initi\'ee
par Berline-Vergne. Pour construire $\eta_\phi$, Bott introduit
$\mu:=\Lcal(V)-\nabla_V$ qui est une section du fibr\'e $\so(\T M)$,
et montre la relation fondamentale
\begin{equation}\label{eq-bott-2}
\nabla\mu=c(V)R \ .
\end{equation}
En introduisant la $1$-forme $\lambda=\frac{(V,-)}{\vert V\vert^2}$,
Bott d\'efinit $\eta_\phi$ sur $M-M^V$ de la mani\`ere suivante
\begin{equation}\label{eq-bott-3}
\eta_\phi:=-\lambda \sum_{k=0}^{N-1} 
C_N^k\,\phi(\underbrace{R,\ldots,R}_{k\ {\rm fois}},\mu,\ldots,\mu)
\,(d\lambda)^{n-k-1}
\end{equation}
o\`u $N$ est le degr\'e du polyn\^ome $\phi$ (suppos\'e inf\'erieur
\`a $n$).

\medskip
 
Montrons maintenant comment la preuve de (\ref{eq-bott-1}) devient
naturelle dans le contexte \'equivariant. On consid\`ere la
sous-alg\`ebre $\Acal(M)^V$, des formes diff\'erentielles invariante
pour l'action infinit\'esimale de $V$, qui est munie de la
d\'erivation $d_V=d-c(V)$. La forme $\phi(R)\in\Acal(M)^V$ n'est 
pas $d_V$-ferm\'ee. Par contre (\ref{eq-bott-2}) montre que la forme 
$\phi(R+\mu)$ est $d_V$-ferm\'ee. Ensuite on remarque que la forme 
$\eta_\phi$ d\'efinie par (\ref{eq-bott-3}) est la composante homog\`ene 
de degr\'e $2n-1$ de la forme
$$
\phi(R+\mu)\frac{\lambda}{d\lambda-1}\ ,
$$
o\`u le d\'enominateur $d\lambda-1$ n'est autre que $d_V\lambda$.
En utilisant la relation fondamentale
$d_V\left(\frac{\lambda}{d_V\lambda}\right)=1$, remarqu\'ee en premier
lieu par Berline-Vergne, on voit que
$$
d_V\left(\phi(R+\mu)\frac{\lambda}{d_V\lambda}\right)=\phi(R+\mu)\quad
{\rm sur}\quad M-M^V\ .
$$

En prenant les composantes homog\`enes de degr\'e maximal $2n$ dans cette
relation on obtient $d\eta_\phi=[\phi(R+\mu)]^{\rm max}$. Sachant que
pour tout polyn\^ome $\phi$ de degr\'e {\em inf\'erieur o\`u \'egal} \`a $n$,
$[\phi(R)]^{\rm max}= [\phi(R+\mu)]^{\rm max}$, on obtient
$d\eta_\phi=[\phi(R)]^{\rm max}$ sur $M-M^V$. $\Box$

\subsection{Localisations de Berline-Vergne et Atiyah-Bott}

Nous avons vu \`a la section pr\'ec\'edente un exemple de localisation
\`a la Berline-Vergne.  Revenons maintenant au cas d'une vari\'et\'e
compacte orient\'ee $M$ munie d'une action d'un groupe de Lie compact
$K$. Soit $\eta$ une forme \'equivariante ferm\'ee \`a coefficients
polyn\^omiaux (on peut prendre indiff\'eremment $\f$).

N. Berline et M. Vergne localisent l'int\'egrale $\int_M\alpha(X)$ sur
la sous-vari\'et\'e $M^X$ des z\'eros du champ de vecteurs $X_M$ 
\cite{B-Vergne1}. On travaille \`a $X\in\kgot$ fix\'e et on proc\`ede 
comme ci-dessus. La forme $\eta(X)$ est ferm\'ee pour la d\'erivation 
$d_V=d-c(X_M)$. Au moyen d'une structure riemannienne $K$-invariante, 
on d\'efinit la $1$-forme $\lambda=(X_M,-)$ qui satisfait la relation
$$
d_V\left(\frac{\lambda}{d_V\lambda}\right)=1 \quad {\rm sur}\quad
M-M^X\ ,
$$
et on localise sur $M^X$ au moyen de $\lambda$.

Sur la sous-vari\'et\'e $M^X$ munie du fibr\'e vectoriel normal
$\Ncal$, on a une action infinit\'esimale $\Lcal(X)$. On munit $\Ncal$
d'une m\'etrique euclidienne $\Lcal(X)$ invariante, et d'une connexion
euclidienne de courbure $R_\Ncal$. L'\'el\'ement $\Lcal(X)+R_\Ncal$
admet un Pfaffien $\det^{1/2}(\frac{\Lcal(X)+R_\Ncal}{-2\pi})\in\Acal(M^X)$ qui est
inversible dans $\Acal(M^X)$ car sa composante homog\`ene de degr\'e
$0$ est constante, \'egale \`a $\det^{1/2}(\frac{\Lcal(X)}{-2\pi})\neq 0$. On a le

\begin{theo}[Berline-Vergne]\label{th-berline-vergne}
\begin{equation}\label{eq-berline-vergne}
\int_M\eta(X)=\int_{M^X}\frac{\eta(X)\vert_{M^X}}
{\det^{1/2}\left(\frac{\Lcal(X)+R_\Ncal}{-2\pi}\right)}\ .
\end{equation}
\end{theo}

\bigskip

{\em Le cas ab\'elien}: Dans le cas o\`u le groupe $K=T$ est un tore,
les sous-vari\'et\'es $M^X$ co\"\i ncident, pour $X$ g\'en\'erique, avec
la sous-vari\'et\'e $M^T$ des points fixes de l'action de $T$. On
d\'efinit la forme d'Euler \'equivariante
$\Eul(\Ncal)(X):=\det^{1/2}(\frac{\Lcal(X)+R_\Ncal}{-2\pi})$ qui est une forme
\'equivariante ferm\'ee sur $M^T$. Le th\'eor\`eme
\ref{th-berline-vergne} donne dans cette situation
\begin{equation}\label{eq-BV-tore}
\int_M\eta(X)=\int_{M^T}\frac{\eta(X)\vert_{M^T}}{\Eul(\Ncal)(X)}\ .
\end{equation}
Cette \'egalit\'e peut \^etre comprise comme une \'egalit\'e de
fonctions sur l'ouvert \break $\{X\in\tgot, M^X=M^T\}$, ou bien comme une
\'egalit\'e alg\'ebrique dans le corps des fractions $\Rcal$ de
$S(\tgot^*)$.

\bigskip

{\em Localisation d'Atiyah-Bott}: Dans le cas ab\'elien, M. Atiyah et
R. Bott obtiennent une localisation directement au niveau de la
cohomologie \cite{Atiyah-Bott}. Pour cela ils montrent que la restriction
$i^*:\Hcal_{T}^{*}(M)\to \Hcal_{T}^{*}(M^T)$ devient un isomorphisme
si l'on ignore la torsion: $\Hcal_{T}^{*}(M))\otimes_{S(\tgot^*)}\Rcal
\stackrel{\sim}{\longrightarrow}\Hcal_{T}^{*}(M^T)\otimes_{S(\tgot^*)}\Rcal $.

Cela entraine que le morphisme
$I_{\rm alg}:\Hcal_{T}^{*}(M)\to\Hcal_{T}^{*}(M)\otimes_{S(\tgot^*)}\Rcal$
d'extension des coefficients se factorise de la mani\`ere suivante

\begin{equation}\label{formuleA-B}
\xymatrix{
\Hcal_{T}^{*}(M)\ar[dr]_{I_{\rm alg}}\ar[r]^-{\Lambda} &
\Hcal_{T}^{*}(M^T)\otimes_{S(\tgot^*)}\Rcal \ar[d]^{i_{*}}\\
 & \Hcal_{T}^{*}(M)\otimes_{S(\tgot^*)}\Rcal
   }
\end{equation}
o\`u le morphisme $\Lambda$ est d\'efini par l'\'equation
$$
\Lambda(\eta)=i^{*}(\eta)\Eul(\Ncal)^{-1},
$$
pour toute classe $\eta\in\Hcal_{T}^{*}(M)$. 

\medskip

L'un des premiers r\'esultats de ma th\`ese a \'et\'e la construction d'un inverse
$\Eul{}^{-1}_\beta(\Ncal)$ de la classe d'Euler \'equivariante $\Eul(\Ncal)$ dans 
l'espace de cohomologie $\Hcal_T^{-\infty}(M^T)$. Ici $\beta\in\tgot$ est tel 
que $M^\beta=M^T$, et
$$
\Eul{}^{-1}_\beta(\Ncal)(X)=\lim_{s\to 0^+}\frac{1}{\Eul(\Ncal)(X+is\beta)}\ .
$$

\subsection{Localisations de Bismut et Witten}

On a vu dans les formules pr\'ec\'edentes que la $1$-forme $\lambda=(\beta_M,-)$
associ\'ee au champ de vecteurs $\beta_M$ joue un grand r\^ole. Dans le cas d'un groupe
compact $K$, on consid\`ere un champ de vecteurs $K$-invariant $\Hcal$ et 
la $1$-forme $\lambda=(\Hcal,-)$ qui est $K$-invariante. Le cas pr\'ec\'edent
correspond \`a $\Hcal=\beta_M$ (invariant seulement si $\beta$ est central).

Donnons dans cette section des r\'esultats valables pour une $1$-forme 
$K$-invariante $\lambda$  quelconque. Cette forme d\'etermine une application 
\'equivariante $\Phi_\lambda:M\to\kgot^*$ satisfaisant la relation
$D\lambda(X)=d\lambda - \langle\Phi_\lambda,X\rangle$.

Si on reprend la m\'ethode de localisation de Berline-Vergne, on cherche \`a
inverser la forme $D\lambda(X)=d\lambda-\langle\Phi_\lambda,X\rangle$. La 
question ici est d'inverser sa composante de degr\'e $0$, 
$X\mapsto\langle\Phi_\lambda,X\rangle$, vue comme un \'el\'ement de 
$S(\kgot^*)\otimes\f(M)$. Cet inverse n'existe pas dans 
un espace de la forme $\Rcal\otimes\f(\Ucal)$, o\`u $\Rcal$ est le corps de
fraction et $\Ucal$ est un ouvert invariant de $M$. Sauf en dimension $1$, car
dans ce cas $\langle\Phi_\lambda,X\rangle=\langle E^*,X\rangle
\langle\Phi_\lambda,E\rangle$ admet pour inverse dans $\Rcal\otimes\f(\Ucal)$,
o\`u $\Ucal=\{\langle\Phi_\lambda,E\rangle\neq 0\}$.

\medskip

Bismut propose une autre m\'ethode de localisation \cite{Bismut}.
\begin{prop}[Bismut]
  Soit $\lambda$ une $1$-forme $K$-invariante sur $M$. Pour toute
  forme $K$-\'equivariante ferm\'ee $\eta(X)$, on a
  $$
  \int_M\eta(X)=\int_M\eta(X)e^{zD\lambda(X)}\quad {\rm pour \ tout
  }\quad X\in\kgot,\ z\in\C\ .
  $$
\end{prop} 

Bismut utilise cette m\'ethode avec $z=t\in \R$ en utilisant la
$1$-forme $\lambda=(\beta_M,-)$, tandis
que Witten \cite{Witten} reprend cette technique avec $z=-i\, t$ pour $t\in\R$.
Expliquons comment ces deux choix donnent des localisations tr\'es
diff\'erentes.

\medskip

{\em Localisation de Bismut}: On fixe $\beta\in\kgot$. Soient $K_\beta$ le sous-groupe
stabilisateur de $\beta$, et $M^\beta$ la sous-vari\'et\'e de $M$ compos\'ee des points 
o\`u $\beta_M$ s'annule. On localise avec la $1$-forme $\lambda=(\beta_M,-)$ qui 
est invariante par rapport \`a $K_\beta$. \`A $X\in\kgot_\beta$ fix\'e, la forme
diff\'erentielle $e^{t\, D\lambda(X)_m},\, m\in M$ est \'egale au
produit $e^{t\, d\lambda_m} \, e^{-t\,(\beta_M,X_M)_m}$. La localisation 
repose sur la d\'ecroissance exponentielle, lorsque $t\to +\infty$, de 
$e^{-t\,(\beta_M,X_M)_m}$ d\`es que $(\beta_M,X_M)_m>0$. Gr\^ace \`a ce 
proc\'ed\'e, Bismut \cite{Bismut} obtient une extension de
(\ref{eq-berline-vergne}) au cadre non-ab\'elien:
\begin{equation}\label{eq-bismut}
\int_M\eta(X)=\int_{M^\beta}\frac{\eta(X)\vert_{M^\beta}}{\Eul(\Ncal)(X)}
\end{equation}
pour $X$ dans un voisinage (assez petit) de $\beta$ dans $\kgot_\beta$.

\bigskip

{\em Localisation de Witten}: Le choix de $z=-i\, t$ permet de
localiser avec plus de souplesse, avec une $1$-forme invariante quelconque.
 La forme \'equivariante $e^{-i\,t\,D\lambda(X)}=e^{-i\,t\,d\lambda}
e^{i\,t\,\langle\Phi_\lambda,X\rangle}$ est le produit d'un terme
polyn\^omial en $t$ avec
$\delta_t:=e^{i\,t\,\langle\Phi_\lambda,X\rangle}\in\f(\kgot,\f(M))$.
Si on se restreint \`a l'ouvert $M-\{\Phi_\lambda=0\}$, on remarque
que $\delta_t$ tend exponentiellement vers $0$ dans l'espace des
fonctions g\'en\'eralis\'ees sur $\kgot$ \`a coefficients dans
$\f(M-\{\Phi_\lambda=0\})$. On montre ainsi le
\begin{lem}[Witten]\label{lem:loc-witten}
  Soit $\chi\in\f(M)$ une fonction test, $0\leq\chi\leq 1$, \'egale
  \`a $1$ au voisinage de $\{\Phi_\lambda=0\}$. Alors pour toute forme
  equivariante ferm\'ee $\eta$, on a
  $\int_{M}(1-\chi)\eta(X)e^{-i\,t\,D\lambda(X)}\longrightarrow 0$
  lorsque $t\to + \infty$. Ainsi
  $$
  \int_{M}\eta(X)=\lim_{t\to+\infty}\int_{M}\chi\,\eta(X)e^{-i\,t\,D\lambda(X)}
  \ .
  $$
  Ici les convergences sont prises dans l'espace des fonctions
  g\'en\'eralis\'ees sur $\kgot$.
\end{lem} 

\medskip

On voit ici les avantages de cette m\'ethode de localisation:
\begin{enumerate}
\item[i)]elle est globale par rapport au param\`etre $X\in\kgot$
  
\item[ii)] on peut prendre des formes $\eta(X)$ \`a coefficients
  $\f$
  
\item[iii)] elle s'effectue aussi bien dans le cas non-ab\'elien
  
\item[iv)] on a une certaine souplesse par rapport aux choix de
  $\lambda$.
\end{enumerate}

N\'eanmoins, cette m\'ethode comporte une s\'erieuse limitation:
l'ensemble $\{\Phi_\lambda=0\}$ n'est g\'en\'eralement pas une
sous-vari\'et\'e.

\subsection{R\'ealisation du programme de Witten}\label{subsec:programme-witten}

Witten a explicit\'e le proc\'ed\'e ci-dessus dans le cadre
hamiltonien \cite{Witten}, en se limitant au calcul au niveau d'une composante lisse
de $\{\Phi_\lambda=0\}$. On reparle de cela dans la prochaine section.

Je d\'ecris ici mes contributions pour d\'evelopper ce type de
localisation \cite{pepthese,pep1,pep2}. Les r\'esultats que j'ai
obtenus sont au carrefour entre la m\'ethode de Berline-Vergne et le
r\'esultat cohomologique de Atiyah-Bott.

\medskip

Je rappelle bri\`evement les notations. Soit $M$ une vari\'et\'e (non
n\'ecessairement compacte) munie d'une action d'un groupe de Lie
compact connexe $K$. On consid\`ere une $1$-forme $\lambda$ sur $M$,
$K$-invariante, et on note $\Phi_\lambda:M\to\kgot^*$ l'application
d\'efinie par $\langle\Phi_\lambda(m),X\rangle=\lambda(X_M)_m$. Le
premier r\'esultat concerne l'inversibilit\'e de la forme
\'equivariante $D\lambda$, o\`u $D\lambda(X)=d\lambda-\langle\Phi_\lambda,X\rangle$.

\begin{lem}[\cite{pep1}] 
  La limite $\lim_{a\to\infty}i\int_0^{a}e^{-i\,t\,D\lambda}dt$
  d\'efinit une forme \'equivariante ferm\'ee \`a coefficients
  g\'en\'eralis\'es sur $M-\{\Phi_\lambda=0\}$. Cette forme
  \'equivariante v\'erifie
  $D\lambda(i\int_0^{\infty}e^{-i\,t\,D\lambda}dt)=1$ sur
  $M-\{\Phi_\lambda=0\}$: on la note $[D\lambda]^{-1}$.
\end{lem}

\medskip

En se pla\c cant dans la cadre des formes \'equivariantes \`a
coefficients g\'en\'eralis\'es, on obtient donc
\begin{equation}\label{eq-fondamentale}
D([D\lambda]^{-1}\lambda)=1 \quad {\rm sur}\quad M-\{\Phi_\lambda=0\} .
\end{equation} 
Une id\'ee naturelle est d'\'etendre (\ref{eq-fondamentale}) \`a $M$.
Pour cela consid\'erons une fonction test $K$-invariante, $\chi$,
\'egale \`a $1$ au voisinage de $\{\Phi_\lambda=0\}$. La forme
$\delta=(1-\chi)[D\lambda]^{-1}\lambda$ est d\'efinie sur $M$, et en
d\'eveloppant $D(\delta)$ on obtient

\begin{lem}[\cite{pep1}] \label{lem-p-lambda}
  La forme $K$-\'equivariante $\loc_\lambda:=\chi + d\chi[D\lambda]^{-1}\lambda$ est
  d\'efinie sur $M$ et satisfait la relation
\begin{equation}\label{eq-p-lambda}
1_M=\loc_\lambda + D(\delta) \quad {\rm sur}\quad M .
\end{equation}
o\`u $1_M$ est la fonction constante \'egale \`a $1$ sur $M$. Ainsi 
pour tout $\eta\in\Hcal_{K}^{*}(M)$, on a
$$
\int_M\eta=\int_{M}\eta\,\loc_\lambda
$$
o\`u l'\'egalit\'e est prise dans l'espace des fonctions 
g\'en\'eralis\'ees $K$-invariantes sur $\kgot$.
\end{lem}

A mon avis l'\'egalit\'e (\ref{eq-p-lambda}) est le paradigme de la
localisation.  On explicite une forme $\loc_\lambda$ \`a support dans
un voisinage de $\{\Phi_\lambda=0\}$ -aussi petit que l'on d\'esire-
qui repr\'esente $1$ en cohomologie. On peut trouver un expos\'e d\'etaill\'e
de l'utilisation de ces formes $\loc_\lambda$, lorsque $K=S^1$, dans \cite{Vergne-Maroc}.

Dans la pratique on est amen\'e \`a consid\'erer chaque `partie' de $\loc_\lambda$. 
On appelle {\em composante} de $\{\Phi_\lambda=0\}$ toute partie 
{\em ferm\'ee} de $\{\Phi_\lambda=0\}$ telle que $\{\Phi_\lambda=0\}-F$ est 
{\em ferm\'ee}. Si $F$ est une composante $K$-invariante de 
$\{\Phi_\lambda=0\}$, on peut d\'efinir la forme \'equivariante ferm\'ee
\`a coefficients g\'en\'eralis\'es
\begin{equation}\label{eq-loc-F}
 \loc_\lambda^F:=\chi^F + 
d\chi^F[D\lambda]^{-1}\lambda
\end{equation}
au moyen d'une fonction $\chi^F$ qui est $K$-invariante, \'egale \`a
$1$ sur un voisinage de $F$, et telle que
$support(\chi^F)\cap\{\Phi_\lambda=0\}=F$. On montre facilement que la classe de
$\loc_\lambda^F$ dans $\Hcal_{K}^{-\infty}(M)$ ne d\'epend pas du choix de la fonction
$\chi^F$. On pr\'ecise alors le lemme \ref{lem:loc-witten} comme suit

\begin{lem}[\cite{pep1}]\label{lem:loc-F}
  Soient $F$ une {\em composante} de $\{\Phi_\lambda=0\}$ et $\chi^F\in\f(M)$ une fonction 
$K$-invariante, \'egale \`a $1$ sur un voisinage de $F$, et telle que $support(\chi^F)\cap
\{\Phi_\lambda=0\}=F$. Alors pour tout $\eta\in\Hcal_{K}^{*}(M)$, on a 
$$
\lim_{t\to + \infty}\int_{M}\chi^F\eta\,e^{-i\,t\,D\lambda}
=\int_{M}\eta\,\loc_\lambda^F
$$
Ici la convergence est prise dans l'espace des fonctions 
g\'en\'eralis\'ees $K$-invariantes sur $\kgot$.  
\end{lem} 
 
Lorsque l'ensemble $\{\Phi_\lambda=0\}$ se d\'ecompose en une union disjointe de
{\em composantes}, $\{\Phi_\lambda=0\}=\cup_i F_i$, la forme
$\loc_\lambda$ s'exprime comme
\begin{equation}\label{eq-loc-decompose}
\loc_\lambda=\sum_i\loc_\lambda^{F_i}\ .
\end{equation}
On sera souvent amen\'e \`a \'etudier chaque composante
$\loc_\lambda^{F_i}$ s\'epar\'ement.

Consid\'erons le cas d'une composante $F$ lisse. On note $i_F$
l'inclusion de $F$ dans $M$, et on suppose que le fibr\'e normal
correspondant $\Ncal_F$ est {\em orient\'e}. On montre de mani\`ere
\'el\'ementaire dans \cite{pep1} que pour tout $\eta\in
\Hcal_{K}^{*}(M)$
\begin{equation}\label{eq-loclisse-F}
  \loc_\lambda^{F}\eta=(i_F)_*\Big(i_F^*(\eta)\Lambda_F\Big)
\end{equation}
o\`u $\Lambda_F\in \Hcal_{K}^{-\infty}(F)$ est un {\em inverse} de la
classe d'Euler \'equivariante du fibr\'e $\Ncal_F$. Dans la pratique,
pour rendre (\ref{eq-loclisse-F}) utilisable, on cherche \`a avoir une
expression explicite de la forme \'equivariante $\Lambda_F$, en
fonction de $\lambda$. Si $F$ est compact, ou si la forme $\eta$ est
\`a support compact sur $M$, on peut int\'egrer l'expression
(\ref{eq-loclisse-F}) :
\begin{equation}\label{eq-localisation-F}
\int_{M}\loc_\lambda^{F}\eta=\int_F i_F^*(\eta)\Lambda_F\ . 
\end{equation}
Ici l'\'egalit\'e est prise dans l'espace de fonctions
g\'en\'eralis\'ees $K$-invariantes sur $\kgot$.  \medskip

Pour terminer, pla\c cons nous dans le cas id\'eal o\`u
$\Ccal=\{\Phi_\lambda=0\}$ est une sous-vari\'et\'e, munie
d'un fibr\'e normal $\Ncal_\Ccal$ orient\'e.  On constate que
(\ref{eq-loclisse-F}) et (\ref{eq-p-lambda}) donnent une factorisation
du morphisme naturel $I:\Hcal_{K}^{*}(M)\to\Hcal_{K}^{-\infty}(M)$:

\begin{equation}\label{formule-lisse}
\xymatrix{
\Hcal_{K}^{*}(M)\ar[dr]_{I}\ar[r]^-{\Lambda} &
\Hcal_{K}^{-\infty}(\Ccal) \ar[d]^{i_{*}}\\
 & \Hcal_{K}^{-\infty}(M)
   }
\end{equation}
o\`u $\Lambda:\Hcal_{K}^{*}(M)\to\Hcal_{K}^{-\infty}(\Ccal)$ est
d\'efini par l'\'equation: $\Lambda(\eta)=i^{*}(\eta)\,
\Lambda_{\Ccal}$.  Comme tout \`a l'heure, la forme $\Lambda_\Ccal\in
\Hcal_{K}^{-\infty}(\Ccal)$ est un inverse de la classe d'Euler
\'equivariante du fibr\'e normal $\Ncal_\Ccal$. Notre 
factorisation (\ref{formule-lisse}) ne requiert pas la compacit\'e 
de la vari\'et\'e $M$. On note les fortes analogies entre 
les factorisation (\ref{formule-lisse}) et  (\ref{formuleA-B}).
 Dans le cadre analytique o\`u je travaille, l'espace de cohomologie
$\Hcal_{K}^{-\infty}(M)$ est l'analogue de
$\Hcal_{T}^{*}(M)\otimes_{S(\tgot^*)}\Rcal$. Remarquons que 
dans certain cas, les morphismes d'extensions des coefficients 
$I_{\rm alg}:\Hcal_{T}^{*}(M)\to\Hcal_{T}^*(M)\otimes_{S^*(\tgot)}\Rcal$ et
$I:\Hcal_{K}^{*}(M)\to\Hcal_{K}^{-\infty}(M)$ 
sont {\em injectifs}: c'est le cas par exemple si l'action est hamiltonienne,
et la vari\'et\'e $M$ est compacte.

Comme dans Atiyah-Bott, notre factorisation donne une localisation des
int\'egrales de formes \'equivariantes. Pour tout
$\eta\in\Hcal_{K}^{*}(M)$ \`a support compact sur $M$, on a
l'\'egalit\'e suivante de fonctions g\'en\'eralis\'ees $K$-invariantes sur
$\kgot$ 
\begin{equation}\label{eq-localisation-generale}
\int_{M}\eta= \int_{\Ccal} i^{*}(\eta)\Lambda_\Ccal \ .
\end{equation}

\medskip

J'ai beaucoup \'etudi\'e cette localisation dans le cadre hamiltonien
o\`u l'hypoth\`ese de lissit\'e de $\{\Phi_\lambda=0\}$ est rarement
v\'erifi\'ee. Je parlerai de cela dans la prochaine section. On va
conclure cette section avec l'illustration de ce proc\'ed\'e lorsque
l'on localise sur les points fixes.

\medskip

{\em Localisation sur les points fixes}. Nous fixons $\beta\in\kgot$, et 
nous proc\'edons \`a la localisation avec la $1$-forme $\lambda=(\beta_M,-)$ qui est
invariante par rapport au sous-groupe stabilisateur $K_\beta$. On voit ici que
$\{\Phi_\lambda=0\}$ co\"\i ncide avec la sous vari\'et\'e $M^\beta$ des points o\`u
le champ de vecteurs $\beta_M$ s'annule.
Le fibr\'e normal $\Ncal_\beta$ de $M^\beta$ dans $M$ est naturellement
orient\'e: soit $\Eul(\Ncal_\beta)$ sa classe d'Euler
$K_\beta$-\'equivariante. On montre dans \cite{pepthese,pep1} que 
la classe de cohomologie $\Lambda_{M^{\beta}}\in
\Hcal_{K_\beta}^{-\infty}(M^{\beta})$ satisfaisant (\ref{eq-loclisse-F})
est \'egale \`a la classe d\'efinie par la forme
\begin{equation}\label{eq-def-euler}
\Eul{}_{\beta}^{-1}(\Ncal_\beta):=\lim_{s\to 0^+}\frac{1}{\Eul(\Ncal_\beta)(X+is\beta)}\ .
\end{equation}
La formule int\'egrale (\ref{eq-localisation-generale}) donne donc
dans ce cas
\begin{equation}\label{eq-euler-beta}
\int_{M}\eta(X)= \lim_{s\to 0^+}\int_{M^{\beta}} \frac{\eta(X)\vert_{M^{\beta}}}
{\Eul(\Ncal_\beta)(X+is\beta)}
\end{equation}
pour toute forme $\eta$, $K_\beta$-\'equivariante, \`a support
compact. Cette \'egalit\'e, qui doit \^etre comprises au sens des
fonctions g\'en\'eralis\'ees, am\'eliore l'expression
(\ref{eq-bismut}) de Bismut.

\bigskip

{\em Polarisation des poids}: La limite (\ref{eq-def-euler}) doit
\^etre comprise comme le proc\'ed\'e de polarisation des poids de
Guillemin-Lerman-Sternberg \cite{G-L-S}, effectu\'e directement au
niveau des formes diff\'erentielles. Consid\'erons l'exemple d'un
fibr\'e vectoriel orient\'e $\Ecal\to M$ muni de l'action d'un tore
$T$ tel que $\Ecal^T=M$. Pour tout $\beta\in\tgot$ tel que
$\Ecal^\beta=M$, on peut d\'efinir $\Eul{}_{\beta}^{-1}(\Ecal)$ par
(\ref{eq-def-euler}). On polarise les poids
$\{\pm\alpha_{1},\cdots,\pm\alpha_{p}\}$ de l'action de $T$ sur les
fibres de $\Ecal$ en convenant que $\alpha_{k}^{+}(\beta)>0\ \forall
k=1,\ldots,p$ . Soit $C_{\beta}$ le c\^one convexe de $\tgot^*$ engendr\'e par
les $\alpha_{k}^{+}$. On montre dans \cite{pepthese,pep1}, que la
transform\'ee de Fourier $\Fcal(\Eul{}^{-1}_{\beta}(\Ecal))$ est
support\'ee par le c\^one $C_{\beta}$. Si de plus l'action de $T$ sur
$\Ecal$ est effective, c'est une mesure {\em continue, localement
  polynomiale} de $C_{\beta}$ dans $\Acal(M)$.

\section{Actions hamiltoniennes et r\'eduction symplectique}

Soit $(M,\omega)$ une vari\'et\'e symplectique compacte munie d'une
action d'un groupe de Lie compact $K$. On suppose que la forme
symplectique est $K$-invariante.  L'action est dite hamiltonienne si
il existe une {\em application moment} $\Phi:M\to\kgot^*$:  l'application $\Phi$ 
est \'equivariante et satisfait les relations
\begin{equation}\label{eq-appl-moment}
d\langle\Phi,X\rangle+\omega(X_M,-)=0\quad \forall X\,\in\,\kgot\ .
\end{equation}

Dans ce cas la forme \'equivariante ferm\'ee
$\omega_\kgot(X):=\omega-\langle\Phi,X\rangle$ est un ant\'ec\'edent
de $\omega$ \`a travers l'application $\Hcal^*_K(M)\to\Hcal^*(M),
\eta(X)\mapsto \eta(0)$: c'est la {\em forme symplectique \'equivariante}.

Dans cette section, on fixe un produit scalaire $K$-invariant sur $\kgot^*$, ce qui  
permet entre autre d'identifier $\kgot$ avec $\kgot^*$. Consid\'erons 
la fonction $\|\Phi\|^2:M\to\R$. Witten propose dans \cite{Witten} de localiser
{\em l'int\'egration} des formes \'equivariantes sur l'ensemble $\crit(\|
\Phi\|^2)$ des points critiques de la fonction $\|\Phi\|^2$. Pour cela
il utilise la $1$-forme
\begin{equation}\label{eq-lambda}
  \lambda=(\Hcal,-)
\end{equation}
o\`u $\Hcal$ est le champ de vecteurs hamiltonien de $\| \Phi\|^2$. 
Dans ce contexte, l'application $\Phi_\lambda:M\to\kgot^*$ d\'efinie par 
$D\lambda(X)=d\lambda -\langle\Phi_\lambda,X\rangle$ v\'erifie
 $\{\Phi_\lambda=0\}=\crit(\|\Phi\|^2)$.  Witten
effectua la localisation sur la composante $\Phi^{-1}(0)$ de
$\crit(\|\Phi\|^2)$ lorsque $0$ est une valeur r\'eguli\`ere de
$\Phi$. Ce r\'esultat a \'et\'e ensuite (re-)d\'emontr\'e par 
Jeffrey et Kirwan \cite{Jeffrey-Kirwan95}. Donnons-en un bref aper\c cu.

Si $0$ est une valeur r\'eguli\`ere de $\Phi$, on peut consid\'erer la
{\em r\'eduction symplectique} en $0$, $\Mcal_0:=\Phi^{-1}(0)/K$, qui
est une\footnote{``Orbifold'' en anglais.} 
V-vari\'et\'e symplectique. Nous avons dans ce contexte deux
morphismes. Le morphisme de Kirwan $\kir_0:\Hcal_K^*(M)\to\Hcal^*(\Mcal_0)$ est
le compos\'e du morphisme de restriction  $\Hcal_K^*(M)\to\Hcal^*_K(\Phi^{-1}(0))$ 
avec l'isomorphisme de Chern-Weil $\Hcal_K^*(\Phi^{-1}(0))
\stackrel{\sim}{\longrightarrow}\Hcal^*(\Mcal_0)$.  
Nous avons aussi l'isomorphisme de Kumar-Vergne
$\kv:\Hcal^*(\Mcal_0)\stackrel{\sim}{\longrightarrow}
\Hcal^{-\infty}_K(\Phi^{-1}(0))$  d\'efinie \`a la fin de 
la section \ref{subsec:modele-cartan} (voir \ref{eq-kv}). Soit $\eta(X)$
une forme $K$-\'equivariante ferm\'ee sur $M$ {\em \`a coefficients
  polynomiaux}, et consid\'erons la forme
$\eta(X)e^{i\omega_\kgot(X)}$. La formule de Jeffrey-Kirwan-Witten
assure que
\begin{eqnarray}\label{eq-JKW}
\lefteqn{\Fcal\left(\int_M\eta\,e^{i\omega_\kgot}\right)=}\\
\nonumber & & (2i\pi)^{\dim K}\ 
\Fcal\left(\int_{\Phi^{-1}(0)}
\kv\circ\kir_0(\eta\,e^{i\omega_\kgot})\right)\quad {\rm au\ voisinage\ de\ } 0,
\end{eqnarray}
o\`u $\Fcal$ d\'esigne la transformation de Fourier.  La fonction
g\'en\'eralis\'ee
$\int_{\Phi^{-1}(0)}\kv\circ\kir_0(\eta\,e^{i\omega_\kgot})$ est
support\'ee en $0$, et sa {\em multiplicit\'e\footnote{Voir la ``convention'' 
que j'ai prise \`a la fin de la section \ref{subsec:modele-cartan}.} par rapport 
\`a la masse de Dirac en $0$} est \'egale \`a
$$
\frac{1}{|S|}\int_{\Mcal_0}\kir_0(\eta\,e^{i\omega_\kgot})
$$
o\`u $|S|$ est le cardinal du stabilisateur g\'en\'erique pour
l'action de $K$ sur $\Phi^{-1}(0)$. Ainsi le terme de gauche de
(\ref{eq-JKW}) est un mesure polynomiale $P(\xi)d\xi$ avec $P(0)=
cst\,\int_{\Mcal_0}\kir_0(\eta\,e^{i\omega_\kgot})$.

\bigskip

Ma th\`ese a consist\'e dans la r\'ealisation de la localisation {\em
  compl\`ete} sur $\crit(\| \Phi\|^2)$ lorsque $K$ est ab\'elien.
Dans \cite{pep1}, je mets en place le proc\'ed\'e g\'en\'eral de localisation que 
j'ai expliqu\'e \`a la section \ref{subsec:programme-witten}, au moyen duquel 
j'\'etend le r\'esultat de ma th\`ese en une localisation cohomologique. 
Je compl\`ete le programme de Witten dans
\cite{pep2}, en obtenant des formules d'induction lorsque le groupe
est non-ab\'elien.

Dans cette premiere sous-section, je vais r\'esumer les r\'esultats
obtenus dans le cadre ab\'elien \cite{pepthese,pep1}. Je
terminerai avec les r\'esultats que j'ai obtenus dans le cas d'un
groupe non-ab\'elien.

\subsection{Le cas ab\'elien}

Dans cette partie, nous supposons que $K=T$ est un tore, et on note
$\Phi_T$ l'application moment. On proc\`ede \`a la
localisation au moyen de la $1$-forme $\lambda^T=(\Hcal^T,-)$, o\`u $\Hcal ^T$
est le hamiltonien de la fonction $\|\Phi_T\|^2$. Dans ce cas
$\{\Phi_{\lambda^T}=0\}=\crit(\|\Phi_T\|^2)$ se d\'ecompose sous la
forme
$$
\crit(\| \Phi_T\|^2)=\bigcup_{\beta\in\Bcal^T}
M^{\beta}\cap\Phi_T^{-1}(\beta)
$$
o\`u $\Bcal^T$ est un sous-ensemble fini de $\tgot\simeq\tgot^*$. 
D'apr\`es (\ref{eq-loc-F}) et (\ref{eq-loc-decompose}), la forme
$\loc_{\lambda^T}$ s'\'ecrit
\begin{equation}\label{eq-loc-abelien}
\loc_{\lambda^T}=\sum_\beta\loc_{\beta}^T\ ,
\end{equation}
o\`u chaque forme \'equivariante $\loc_{\beta}^T$ est support\'ee sur
un voisinage (petit) de $M^{\beta}\cap\Phi_T^{-1}(\beta)$. En
g\'en\'eral les sous-ensembles $M^{\beta}\cap\Phi_T^{-1}(\beta)$ ne
sont pas lisses. Pour rem\'edier \`a cela on consid\`ere les
applications $\Phi_T-\esp$ et les $1$-formes $\lambda^T_\esp$
correspondantes.  On aura deux points de vue: l'un {\em local} et
l'autre {\em global}. Du point de vue local, on cherche, en prenant $\esp$
petit, \`a obtenir une `d\'esingularisation' de
$M^{\beta}\cap\Phi_T^{-1}(\beta)$. Un lemme de d\'eformation
\cite{pep1,pep2} donne alors une expression cohomologique de
$\loc_{\beta}^T$ correspondant \`a cette d\'eformation.

Dans l'autre point de vue, on consid\`ere des $\esp$ g\'en\'eraux tel que
$\crit(\| \Phi_T-\esp\|^2)$ est lisse, et on obtient une formule
cohomologique globale du type (\ref{formule-lisse}).

\medskip

Dans ma th\`ese j'obtiens une d\'ecomposition de $\crit(\| \Phi_T-\esp\|^2)$
param\'etr\'ee par une collection $\Acal$ de sous-espaces affines de
$\tgot^*$ d\'efinie de la mani\`ere suivante. La vari\'et\'e $M$ \'etant compacte,
l'action de $T$ sur $M$ poss\`ede un nombre fini de types d'orbites: soient 
$T_1,\ldots,T_r$ les sous-groupes de $T$, stabilisateurs de points de $M$.
Pour chaque $l=1,\ldots,r$ on note $Z^k_l, k=1,\ldots,n_l$ les composantes connexes
de $M^{T_l}$ qui ont pour stabilisateur g\'en\'erique le sous-groupe $T_l$. Les $Z^k_l$
sont des sous-vari\'et\'es symplectiques $T$-invariantes de $M$. 
Le th\'eor\`eme de convexit\'e
d'Atiyah-Guillemin-Sternberg assure que les $\Phi_T(Z^k_l)$ sont des polytopes convexes
de $\tgot^*$. La collection $\Acal$ est l'ensemble des sous-espaces affines
de $\tgot^*$ engendr\'es par les polytopes $\Phi_T(Z^k_l)$. Pour chaque 
$\Delta\in\Acal$, on note $T_\Delta$ le
sous-tore de $T$ d'alg\`ebre de Lie $\tgot_\Delta:=\{X\in\tgot,
\langle a-b,X\rangle=0\ \forall a,b\, \in\Delta\}$.

\begin{prop}[\cite{pepthese}]
  Pour tout $\esp\in\tgot^*$, les points critiques de $\|
  \Phi_T-\esp\|^2$ se mettent sous la forme
\begin{equation}\label{eq-point-critique}
\crit(\| \Phi_T-\esp\|^2)=\bigcup_{\p\in\Acal} M^{T_\Delta}\cap
\Phi^{-1}_T(\beta(\esp,\Delta))
\end{equation}
o\`u $\beta(\esp,\Delta)$ est le projet\'e orthogonal de $\esp$ sur
$\Delta$. Pour $\esp$ g\'en\'erique, l'ensemble $\crit(\|
\Phi_T-\esp\|^2)$ est une sous-vari\'et\'e de $M$, la r\'eunion
(\ref{eq-point-critique}) est disjointe, et le groupe $T/T_\Delta$
agit localement librement sur
$C_\Delta^\esp:=M^{T_\Delta}\cap\Phi^{-1}_T(\beta(\esp,\Delta))$.
\end{prop}

\medskip
 
{\em R\'esulat local}. Consid\'erons une composante critique
$M^{\beta}\cap\Phi^{-1}_T(\beta)$ munie d'un voisinage $\Ucal$ tel que
$\crit(\| \Phi_T\|^2)\cap\Ucal = M^{\beta}\cap\Phi^{-1}_T(\beta)$.
Maintenant effectuons la d\'eformation $\Phi_T\to\Phi_T-\esp$.  Pour
$\Ucal$ convenablement choisi on voit que
$$
\crit(\| \Phi_T - \esp\|^2)\cap\Ucal=\bigcup_{\beta(0,\p)=\beta}
C_\Delta^\esp
$$
pour tout $\esp$ suffisamment petit. Ici l'union est restreinte aux
sous-espaces affines $\p$ tels que la projection orthogonale de $0$
sur $\p$ est \'egale \`a $\beta$. Si $\esp$ est de plus g\'en\'erique,
l'union pr\'ec\'edente est disjointe, et on d\'efinit pour chaque $\p$
la forme $\loc_\p^\esp$, qui est support\'ee au voisinage de
$C_\Delta^\esp$ (voir (\ref{eq-loc-F})).  Gr\^ace \`a un lemme de
d\'eformation (Proposition 2.6 dans \cite{pep2}), on a
\begin{equation}
  \label{eq:loc-deformation}
  \loc^T_\beta=\sum_{\beta(0,\p)=\beta}\loc_\p^\esp \quad\quad {\rm dans}\quad
\Hcal_T^{-\infty}(M).
\end{equation}

\medskip

Nous donnons maintenant l'expression de la localisation obtenue avec
chaque forme $\loc_\p^\esp$. Fixons $\esp$ g\'en\'erique, que l'on ne
suppose plus \^etre petit.  Pour chaque $\Delta$, le quotient
$\Mcal_\Delta^\esp:=C_\Delta^\esp/(T/T_\Delta)$ est une V-vari\'et\'e
munie d'une action triviale du tore $T_{\Delta}$. Dans ce contexte,
nous avons deux morphismes. Le premier est celui de Kirwan
$$
\kir_\p:\Hcal^*_T(M)\to\Hcal^*_{T_\Delta}(\Mcal_\Delta^\esp)
$$
qui est le compos\'e de la restriction
$\Hcal^*_T(M)\to\Hcal^*_{T}(C_\Delta^\esp)$ et de l'isomorphisme de
Chern-Weil
$\Hcal^*_{T}(C_\Delta^\esp)\stackrel{\sim}{\longrightarrow}\Hcal^*_{T_\Delta}
(\Mcal_\Delta^\esp)$. L'autre est l'isomorphisme de Kumar-Vergne (voir
(\ref{eq-kv})).
$$
\kv_\Delta:\Hcal^{-\infty}_{T_\Delta}
(\Mcal_\Delta^\esp)\stackrel{\sim}{\longrightarrow}\Hcal^{-\infty}_{T}(C_\Delta^\esp)\ 
.
$$

Soient $\beta_{\p}:= \beta(\esp,\p)-\esp\in\tgot_{\p}$ et $N_{\p}$ le
fibr\'e normal de $M^{T_{\p}}$ dans $M$ restreint \`a
$C^{\epsilon}_{\p}$. Consid\'erons le $V$-fibr\'e vectoriel
$\Ecal_\p:=N_{\p}/(T/T_\p)$ sur $\Mcal_\Delta^\esp$.  Ces donn\'ees
nous permettent de d\'efinir $\Eul{}^{-1}_{\beta_{\p}}(\Ecal_{\p})\in
\Hcal^{-\infty}_{T_\Delta}(\Mcal_\Delta^\esp)$.

On montre dans \cite{pep2} que la forme $\loc_\p^\esp$ est
d\'etermin\'ee par la relation
\begin{equation}\label{eq:calcul-loc-delta}
  \frac{1}{(2i\pi)^{\dim\p}}\ \loc_\p^\esp \,\eta=(i_\p)_*\left(\kv_\p\Big(
\kir_\p(\eta)\Eul{}^{-1}_{\beta_{\p}}(\Ecal_{\p})\Big)\right)
\end{equation}
pour tout $\eta\in\Hcal^*_{T}(M)$. Ici $i_\p$ d\'esigne l'inclusion de 
$C_\Delta^\esp$ dans $M$. Il faut comprendre cette formule comme un 
{\em m\'elange} du cas o\`u $\p$ est un point, et celui o\`u $\p=\tgot^*$.

\medskip

$\bullet\ \p=\{p\}$ est un sommet du polytope $\Phi_T(M)$. Dans ce cas
$C_\Delta^\esp=\phi^{-1}_T(p)$ est une composante connexe $F$ de $M^T$,
et le tore $T_\p$ est \'egal \`a $T$. L'expression
(\ref{eq:calcul-loc-delta}) devient
$$
\loc_{\{p\}}^\esp \eta=(i_F)_*\left(i_F^*(\eta)
  \Eul{}^{-1}_{\beta_{p}}(\Ncal_{F})\right)\ ,
$$
avec $\beta_p=p-\esp$.

\medskip

$\bullet\ \p=\tgot^*$. Ici $C_\Delta^\esp=\phi^{-1}_T(\esp)$ et
$\Mcal_\p^\esp$ correspond \`a la vari\'et\'e  r\'eduite
$\Mcal_\esp:=\phi^{-1}_T(\esp)/T$.  L'expression
(\ref{eq:calcul-loc-delta}) correspond dans ce cas \`a la localisation
de Jeffrey-Kirwan-Witten
\begin{equation}
  \label{eq:loc-en-esp}
 \frac{1}{(2i\pi)^{\dim T}}\  \loc_{\tgot^*}^\esp\, \eta =
(i_\esp)_*\circ\kv_\esp\circ\kir_\esp(\eta)\ .
\end{equation}
Il est int\'eressant de visualiser le morphisme
$(i_\esp)_*\circ\kv_\esp\circ\kir_\esp$:
\begin{equation}
  \label{eq:morphisme-en-esp}
  \Hcal^{*}_T(M)\stackrel{\kir_\esp}{\longrightarrow}
\Hcal^{*}(\Mcal_\esp)\stackrel{\kv_\esp}{\longrightarrow}
\Hcal^{-\infty}_T(\Phi^{-1}_T(\esp))\stackrel{(i_\esp)_*}{\longrightarrow}
\Hcal^{-\infty}_T(M)\ .
\end{equation}

\bigskip

{\em R\'esultat gobal}. Pour $\esp$ g\'en\'erique, l'expression
(\ref{eq:calcul-loc-delta}) obtenue pour chaque $\p$ d\'efinit une
localisation globale au niveau de la cohomologie (voir
\ref{formule-lisse}). On a la factorisation suivante
\begin{equation}\label{formule-lisse-tore}
\xymatrix{
\Hcal_{T}^{*}(M)\ar[dr]_{I}\ar[r]^-{{\bf kir}} &
\oplus_{\p}\Hcal_{T_\p}^{*}(\Mcal^\esp_\p) \ar[d]^{{\bf j}}\\
 & \Hcal_{T}^{-\infty}(M)
   }
\end{equation}
o\`u ${\bf kir}:=\oplus \kir_\Delta$, et le morphisme ${\bf j}$ est
d\'efini par
$$
{\bf j}(\eta)=i_*\left(\sum_\p \,(2i\pi)^{\dim\p}\ 
  \kv_\p\left(\eta\,\Eul{}^{-1}_{\beta_{\p}}(\Ecal_{\p})\right)\right)\ 
.
$$
Rappelons que le morphisme $I$ d'extension des coefficients est injectif, ce qui entraine 
que ${\bf kir}$ est aussi injectif.

 \medskip

 Dans la prochaine sous-section, on va voir que notre {\em r\'esultat
   local} permet d'obtenir assez facilement les formules de saut de V.
 Guillemin et J. Kalkman \cite{cras04}.

 \subsection{Les formules de saut de Guillemin-Kalkman}

 On travaille dans les m\^emes conditions qu'\`a la section
 pr\'ec\'edente. Par commodit\'e, on suppose que l'action de $T$ sur
 $M$ est effective. On se fixe pour cette section une forme $\eta\in
 \Hcal_{T}^{*}(M)$. A chaque valeur $\xi\in\tgot^*$ r\'eguli\`ere de
 $\Phi_T$, on peut associer
 $$
 I(\xi):=\frac{1}{|S_\xi|}\int_{\Mcal_\xi}\kir_{\xi}(\eta)
 $$
 o\`u $|S_\xi|$ est le cardinal du stabilisateur g\'en\'erique de T
 sur $\Phi^{-1}(\xi)$.  On va voir que notre technique de localisation
 permet de retrouver les propri\'et\'es de $\xi\to I(\xi)$: tout
 d'abord le fait que cet application est localement constante, et
 ensuite les formules de saut de Guillemin-Kalkman.

 Tout d'abord int\'egrons (\ref{eq:loc-en-esp}), en prenant $\esp=\xi$.
 La fonction g\'en\'eralis\'ee $\int_M\loc_{\tgot^*}^\xi \eta$ est
 \'egale \`a $(2i\pi)^{\dim T}\int_{\Phi_T^{-1}(\xi)}\kv_\xi\circ\kir_\xi(\eta)$. 
 On a remarqu\'e \`a la fin de la section \ref{subsec:modele-cartan} que
 cette derni\`ere fonction g\'en\'eralis\'ee est support\'ee en $0$, et
 que sa multiplicit\'e par rapport \`a la masse de Dirac en $0$ est
 \'egale \`a $\frac{(2i\pi)^{\dim T}}{|S_\xi|}\int_{\Mcal_\xi}\kir_{\xi}(\eta)$. Nous
 avons donc une caract\'erisation de $I(\xi)$, comme la multiplicit\'e
 par rapport \`a la masse de Dirac en $0$ de la fonction
 g\'en\'eralis\'ee $\frac{1}{(2i\pi)^{\dim T}}\int_M\loc_{\tgot^*}^\xi \eta$.

 \medskip

 {\em L'application $I$ est localement constante.} Consid\'erons une
 valeur r\'eguli\`ere $\xi\in\tgot^*$ de $\Phi$. Quitte \`a modifier
 $\Phi_T$ en $\Phi_T-\xi$, on peut supposer que $\xi=0$. Le seul
 sous-espace affine $\p$ qui contient la valeur r\'eguli\`ere $0$ est
 $\tgot^*$.

 Nous proc\'edons \`a la d\'eformation de $0$ en une valeur proche
 $\xi'$, et utilisons la relation (\ref{eq:loc-deformation}) avec
 $\beta=0$:
 $$
 \loc_{\tgot^*}^0=\loc_{\tgot^*}^{\xi'}\quad\quad {\rm dans}\ \ 
 \Hcal_{T}^{-\infty}(M)\ ,
 $$
 ce qui d\'emontre que $I(0)=I(\xi')$.

 \medskip

 {\em Les formules de saut.} Consid\'erons un hyperplan $\p_o\in\Bcal$
 s\'eparant deux r\'egions de valeurs r\'eguli\`eres de $\Phi$. Dans ce
 cas le sous-tore $T_{\p_o}$ est de dimension $1$. Soit $a\in\p_o$ 
 tel que $C_{\p_o}^a:=\Phi^{-1}(a)\cap M^{T_{\p_o}}$ est
 une sous vari\'et\'e munie d'une action localement libre de
 $T/T_{\p_o}$.  Quitte \`a modifier $\Phi_T$ en $\Phi_T-a$, on peut
 supposer que $a=0$. Les seuls sous-espaces affines $\p$ qui
 contiennent $0$ sont $\p_o$ et $\tgot^*$.

 Nous utilisons la relation (\ref{eq:loc-deformation}) avec $\beta=0$,
 et l'appliquons \`a des valeurs r\'eguli\`eres de $\Phi$, $\xi_\pm$, 
proches de $a(=0)$, qui sont de chaque c\^ot\'e de $\p_o$, et
 qui se proj\`etent toutes deux sur $a$. Ceci donne
 $\loc_{\tgot^*}^{\xi_+}+\loc_{\p_o}^{\xi_+}=\loc_{\tgot^*}^{\xi_-}+\loc_{\p_o}^{\xi_-}$
 dans $\Hcal^{-\infty}_T(M)$. Les relations (\ref{eq:calcul-loc-delta})
 donnent apr\`es int\'egration

 \begin{eqnarray}\label{eq-GK-1}
 \lefteqn{ \int_{\phi^{-1}(\xi_+)}\kv_{\tgot^*}\circ\kir_{\xi_+}(\eta)\ - \ 
 \int_{\phi^{-1}(\xi_-)}\kv_{\tgot^*}\circ\kir_{\xi_-}(\eta)}\\
 \nonumber &=& \frac{1}{2i\pi}\int_{C_{\p_o}^a}  \kv_{\p_o}
 \left(\kir_{\p_o}(\eta)\Big[\Eul{}^{-1}_{\beta}(\Ecal_{\p_o})-
 \Eul{}^{-1}_{-\beta}(\Ecal_{\p_o})\Big]\right)
 \end{eqnarray}
 o\`u $\beta=\xi_+-\xi_-$ est un vecteur non-nul de $\tgot_{\p_o}$.

 \medskip

 Dans notre situation, $\tgot_{\p_o}$ est de dimension $1$. De plus, la
 forme d'Euler \'equivariante $\Eul{}(\Ecal_{\p_o})$, vue comme un
 polyn\^ome sur $\tgot_{\p_o}$ \`a valeur dans $\Acal(\Mcal_{\p_o}^a)$,
 est inversible de mani\`ere $\f$ sur $\tgot_{\p_o}-\{0\}$. Comme
 $\Eul{}^{-1}_{\pm\beta}(\Ecal_{\p_o})$ sont deux inverses de
 $\Eul{}(\Ecal_{\p_o})$, on en conclut que $\Eul{}^{-1}_{\beta}
 (\Ecal_{\p})-\Eul{}^{-1}_{-\beta}(\Ecal_{\p})$ est une fonction
 g\'en\'eralis\'ee support\'ee en $0$. On peut alors d\'efinir une
 application {\em r\'esidu}
 \begin{equation}\label{eq-residu}
 {\rm Res}_{\p_o}:\Hcal^*_{T_{\p_o}}(\Mcal_{\p_o}^a)\to\Hcal^*(\Mcal_{\p_o}^a)
 \end{equation}
 de la mani\`ere suivante. Pour toute forme $\eta$, son r\'esidu 
${\rm   Res}_{\p_o}(\eta)$ est la composante par rapport \`a la masse de Dirac
 en $0$ de la fonction g\'en\'eralis\'ee \break
 $(2i\pi)^{-1}\eta\,(\Eul{}^{-1}_{\beta}(\Ecal_{\p_o})-
\Eul{}^{-1}_{-\beta}(\Ecal_{\p_o}))$.

 En consid\'erant les multiplicit\'es par rapport \`a la masse de Dirac
 en $0$ dans l'\'egalit\'e (\ref{eq-GK-1}), on obtient la formule de
 saut de Guillemin-Kalkman \cite{Guillemin-Kalkman}
 $$
 I(\xi_+)- I(\xi_-) =
 \frac{1}{|S_\p|}\int_{\Mcal_{\p_o}^a} {\rm
   Res}_{\p_o}\left(\kir_{\p_o}(\eta)\right).
 $$

\subsection{Le cas non-ab\'elien}

Dans cette partie on consid\`ere l'action hamiltonienne d'un groupe de
Lie compact $K$.  Soit $T$ un tore maximal de $K$, et $W$ le groupe de
Weyl. On note respectivement $\Phi_K$ et $\Phi_T$ les applications
moments pour les actions de $K$ et $T$ sur $M$. Rappelons que 
$\Phi_T$ est le compos\'e de $\Phi_K$ avec la projection $\kgot^*\to
\tgot^*$.

Dans le contexte non-ab\'elien, on ne peut plus esp\'erer
`d\'esingulariser' les points critiques de $\|\Phi_K\|^2$, en prenant
$\Phi_K-\esp$, car cette derni\`ere application n'est g\'en\'eralement
plus $K$-\'equivariante. Dans \cite{pep2}, on a contourn\'e cette
difficult\'e en d\'emontrant une propri\'et\'e d'{\em induction}.

Les points critiques de $\|\Phi_K\|^2$ admettent la d\'ecomposition
disjointe \cite{Kirwan.84}
\begin{equation}\label{eq-points-critiques}
\crit(\|\Phi_K\|^2)=\bigcup_{\beta\in\Bcal^K}K(M^{\beta}\cap\Phi_K^{-1}(\beta))  
\end{equation}
o\`u $\Bcal^K$ est une partie finie d'une chambre de Weyl $\tgot^*_+$. Pour les
points critiques de $\|\Phi_T\|^2$, on a de m\^eme
$\crit(\|\Phi_T\|^2)=\bigcup_{\beta\in\Bcal^T}M^{\beta}\cap\Phi_T^{-1}(\beta)$.
Cette deni\`ere d\'ecomposition est \'equivariante par rapport au
groupe de Weyl.  Ainsi $\Bcal^T=W.(\Bcal^T\cap \tgot^*_+)$, et
$\Bcal^K\subset\Bcal^T\cap\tgot^*_+$.

\medskip

Pour chaque $\beta\in\Bcal^K$, on d\'efinit la forme \'equivariante
ferm\'ee $\loc_{\beta}^{K}$ d\'efinie avec la $1$-forme
(\ref{eq-lambda}), et support\'ee sur voisinage de
$K(M^{\beta}\cap\Phi_K^{-1}(\beta))$ (voir (\ref{eq-loc-F})). Ici le
lemme (\ref{lem-p-lambda}) donne la partition
\begin{equation}
  \label{eq:partition}
  1_M=\sum_{\beta\in\Bcal^K}\loc_{\beta}^{K}\quad\quad {\rm
  dans}\quad\Hcal_K^{-\infty}(M)\ .
\end{equation}
On proc\`ede de m\^eme avec le tore $T$. Pour chaque $\beta\in\Bcal^T$,
on d\'efinit $\loc_{\beta}^{T}\in\Hcal_{T}^{-\infty}(M)$ support\'ee sur 
un voisinage de $M^{\beta}\cap\Phi_T^{-1}(\beta)$.

\medskip

Soit $\indt :\Hcal^{-\infty}_{T}(M) \longrightarrow \Hcal^{-\infty}_{K}(M)$ 
le morphisme d'induction d\'efini par Kumar-Vergne dans \cite{Kumar-Vergne}. 
Lorsque $M=\{.\}$, c'est le morphisme d'induction usuel
$\indt :\fgene(\tgot) \longrightarrow \fgene(\kgot)^K$. Dans la suite,
le polyn\^ome $Y\in\tgot\mapsto \det_{\kgot/\tgot}(ad(Y))$ est not\'e
$\Pi_{\kgot/\tgot}^2$. La formule int\'egrale de Weyl s'exprime au
moyen de $\indt$ sous la forme suivante: $f=|W|^{-1}\indt(
\Pi_{\kgot/\tgot}^2 f|_\tgot)$ pour tout $f\in\f(\kgot)^K$.

Ce morphisme d'induction est `fonctoriel' de la mani\`ere suivante. Si
$M$ est orient\'ee, le diagramme suivant est commutatif.
\begin{equation}\label{ind-fonctoriel-point}
\xymatrix@R=15mm@C=2cm
{
\Hcal_{T}^{-\infty}(M) \ar[d]_{\int_{M}}    \ar[r]^{\indt}  &  
\Hcal_{K}^{-\infty}(M) \ar[d]_{\int_{M}} \\
\fgene(\tgot) \ar[r]^{\indt}     & \fgene(\kgot)^K.   \\ 
}
\end{equation}

L'un des r\'esultat principaux de \cite{pep2} est la relation
d'induction suivante. Pour chaque $\beta\in\Bcal^K$, on note $W_\beta$
le sous-groupe des \'el\'ements de $W$ laissant $\beta$ fixe, et
$|W_\beta|$ son cardinal.

\begin{prop}\cite{pep2} \label{theo-induction}
  Les formes \'equivariantes 
  $(\loc_{\beta}^{T})_{\beta\in\Bcal^T}$, et
  $(\loc_{\beta}^{K})_{\beta\in\Bcal^K}$ satisfont les relations
  suivantes: pour tout $\beta\in\Bcal^K$ on a
  $$
  \loc_{\beta}^K =\frac{1}{|W_{\beta}|}
  \indt\Big(\Pi_{\kgot/\tgot}^2\,\loc_{\beta}^T\Big) \quad {\rm dans}
  \quad \Hcal^{-\infty}_{K}(M)\ .
  $$
  Sinon, pour $\beta\in\Bcal^T$ tel que $\beta\notin W.\Bcal^K$, on
  a $\indt\left(\Pi_{\kgot/\tgot}^2\,\loc_{\beta}^T\right)= 0$ dans
  $\Hcal^{-\infty}_{K}(M)$.
\end{prop}

Le cas $\beta=0$ est particuli\`erement int\'eressant. L'\'egalit\'e
$\indt(\Pi_{\kgot/\tgot}^2\,\loc_{0}^T)= |W|\loc_{0}^K$ montre que
pour tout $\eta\in\Hcal^*_K(M)$, on a l'\'egalit\'e de fonctions
g\'en\'eralis\'ees suivante
$$
|W|\int_M\loc_0^K\eta=\indt\left(\Pi_{\kgot/\tgot}^2\int_M\loc_0^T\eta|_\tgot\right)\ 
.
$$
Supposons que $0$ est une valeur r\'eguli\`ere de $\Phi_K$ et de
$\Phi_T$. Alors les deux termes de l'\'egalit\'e pr\'ec\'edente sont
des fonctions g\'en\'eralis\'ees support\'ees en $0$. Si on prend les
multiplicit\'es de la masse de Dirac en $0$, on obtient
\begin{equation}\label{eq-martin}
\frac{|W|}{|S^K|}\int_{\Mcal_0^K}\kir_{0}^K(\eta)=
\frac{1}{|S^T|}\int_{\Mcal_0^T}\kir_{0}^T(\Pi_{\kgot/\tgot}^2\eta|_\tgot)
\end{equation}

Cette derni\`ere relation a aussi \'et\'e obtenue par Martin
\cite{Martin} avec des techniques diff\'erentes.

\bigskip

Nous terminons cette section en donnant une application de la
proposition \ref{theo-induction}.

\medskip

{\em Les fonctions de partition de Witten}. Dans \cite{Witten}, Witten
introduit les fonctions de partition $Z$ d\'efinies par
$$
Z(u):=\int_{\kgot}e^{-u\|X\|^2/2}\left(\int_M\eta
  e^{i\omega_\kgot}\right)(X)dX\ , \quad u>0,
$$
o\`u $\eta$ est une forme $K$-\'equivariante ferm\'ee \`a
coefficients {\em polyn\^omiaux}, et $\omega_\kgot$ est la forme
symplectique \'equivariante.

En supposant que $0$ est une valeur r\'eguli\`ere de $\Phi_K$, Witten
donne les premiers termes du d\'eveloppement asymptotique de $Z(u)$
lorsque $u\to 0^+$ : $Z(u)=Z_0(u) + \Ocal(e^{-\rho/u})$ avec $\rho>0$
et $Z_0(u)$ est un polyn\^ome. Il obtient une expression explicite de
$Z_0(u)$ sous la forme suivante
$$
Z_0(u)=cst\int_{\Mcal_0^K}\kir_0(\eta
e^{i\omega_\kgot})e^{-u\|\Omega\|^2/2}.
$$
Dans cette expression $\Omega$ est la courbure du fibr\'e
principal $\Phi^{-1}_K(0)\to\Mcal_0^K$. Jeffrey et Kirwan ont
donn\'e une d\'emonstration (rigoureuse) de cette formule dans
\cite{Jeffrey-Kirwan95}. Pour le reste $\Ocal(e^{-\rho/u})$, Witten
conjectura une expression de la forme $\sum_F Z_F(u)$, param\'etr\'ee
par les composantes de $\crit(\|\Phi_K\|^2)$.

Dans \cite{pep2}, j'obtiens le d\'eveloppement asymptotique complet de
$Z(u)$, {\em sans supposer que $0$ est une valeur r\'eguli\`ere de
$\Phi_K$}. Rappelons qu'un ensemble fini $\Bcal^K\subset\tgot^*$
param\`etre $\crit(\|\Phi_K\|^2)$.

\begin{prop}[\cite{pep2}] La fonction $Z(u)$ admet la d\'ecomposition suivante
  $$
  Z(u)=u^{-N}\sum_{\beta\in\Bcal^K}e^{-\frac{\|\beta\|^2}{2u}}h_\beta(\sqrt{u}),\quad
 u>0 .
  $$
  Dans cette expression les fonctions $h_\beta:\R\to\C$ sont $\f$,et la fonction $h_0$
  est {\em toujours polynomiale}.  Si $0$ est une valeur
  r\'eguli\`ere de $\Phi_K$, la fonction $Z_0(u)=u^{-N} h_0(\sqrt{u})$
  est polynomiale. L'entier $N$ est \'egal \`a $\dim K + d(\eta)$,
  o\`u $d(\eta)$ est le degr\'e (en X) de $\eta$.
\end{prop}

\subsection{Mesures de Duistermaat-Heckman et orbites coadjointes}

Nous travaillons toujours dans le contexte d'une action hamiltonienne
d'un groupe de Lie compact $K$ sur une vari\'et\'e symplectique
$(M,\omega)$ de dimension $2n$.  Nous supposons ici que l'application
moment $\Phi$ est propre. Un invariant int\'eressant est l'image
directe $\mdh(M):= \Phi_*(\frac{\omega^n}{n!})$ de la mesure de
Liouville par l'application moment: c'est une mesure $K$-invariante
sur $\kgot^*$ support\'ee par $\Phi(M)$. Lorsque $M$ est compacte,
$\mdh(M)$ s'exprime comme la transform\'ee de Fourier
\begin{equation}
  \label{eq:mdh-compact}
  \mdh(M)=\frac{1}{i^n}\Fcal\left(\int_M e^{i\omega_\kgot}\right)\ ,
\end{equation}
o\`u $\omega_\kgot$ est l'application moment \'equivariante.

\bigskip

{\em Cas ab\'elien}.

\medskip

Dans \cite{D-H}, Duistermaat et Heckman montrent que pour $K=T$
ab\'elien et $M$ compacte, la mesure $\mdh(M)$ est polynomiale sur
chaque sous-polytope de $\Phi(M)$ qui est compos\'e de valeurs
r\'eguli\`eres de $\Phi$.  Si on utilise la formule de Berline-Vergne
dans (\ref{eq:mdh-compact}), et un vecteur $\beta\in\tgot$ tel que
$M^\beta= M^T$, on obtient
\begin{equation}
  \label{eq:mdh-decomp}
\mdh(M)=\sum_{F\subset M^T} \, \delta_{\Phi(F)} * \mdh_\beta^F\,
\end{equation}
Dans cette expression on somme sur les composantes connexes de $M^T$,
$*$ est le produit de convolution et $\delta_{\Phi(F)}$ est la mesure
de Dirac en $\Phi(F)$. Chaque terme $\mdh_\beta^F$ est une mesure
localement polynomiale support\'ee par le c\^one $\sum \R\alpha_i^+$,
o\`u $\alpha_i^+$ sont les poids polaris\'es par $\beta$, de l'action
de $T$ sur $\T M|_F$.

Prato et Wu ont montr\'e que l'expression (\ref{eq:mdh-decomp}) est
encore valable lorsque $M$ est {\em non-compacte} si les deux
conditions suivantes sont satisfaites:
\begin{enumerate}
\item la fonction $\langle\Phi,\beta\rangle$ est propre et born\'ee inf\'erieurement\\
\item l'ensemble des points fixes $M^T$ est fini.
\end{enumerate}
Dans \cite{pep2}, je montre que la deuxi\`eme condition n'est pas
n\'ecessaire, c'est \`a dire que d'une part on peut autoriser des
sous-vari\'et\'es de points fixes, et d'autre part l'ensemble $M^T$
peut avoir un nombre de composantes connexes infini.

\bigskip

{\em Un cas non-ab\'elien: les orbites coadjointes de groupes
 r\'eels semi-simples}.

\medskip

Un exemple important d'action hamiltonienne sur des vari\'et\'es
non-compactes, mais poss\'edant une application moment {\em propre}
est le suivant. Soit $G$ un groupe r\'eel semi-simple connexe
d'alg\`ebre de Lie $\ggot$. Ici $K$ d\'esigne un sous-groupe compact
maximal de $G$. Soit $\Ocal\subset\ggot^*$ une orbite coadjointe de
$G$, munie de la forme symplectique de Kirillov-Kostant-Souriau. On
note $\Phi_\Ocal:\Ocal\to\kgot^*$ l'application moment relative \`a
l'action de $K$: c'est le compos\'e de l'inclusion $\Ocal\croc\ggot^*$ 
avec la projection $\ggot^*\to\kgot^*$.

On remarque tout d'abord que l'application moment $\Phi_\Ocal$ est
{\em propre} si et seulement si $\Ocal$ est {\em ferm\'ee} dans
$\ggot^*$ \cite{pep3}. Dans ce cas l'\'egalit\'e (\ref{eq:mdh-compact}) est encore
valable, parce que l'int\'egrale
$$
F_\Ocal(X)=\int_\Ocal e^{i\omega_\kgot(X)}
$$
d\'efinit une fonction g\'en\'eralis\'ee {\em temp\'er\'ee}. La
fonction g\'en\'eralis\'ee $F_\Ocal$ a \'et\'e calcul\'ee dans
plusieurs situations:

\begin{enumerate}
\item[i)] Lorsque $\Ocal$ est elliptique r\'eguli\`ere, Rossmann
  \cite{Ross} a calcul\'e $F_\Ocal(X)$ pour tout \'el\'ement
  r\'egulier de $\kgot$
  (i.e. l\`a o\`u $F_\Ocal$ est lisse).\\
  
\item[ii)] Lorsque $\Ocal$ est r\'eguli\`ere (i.e. de dimension
  maximale), $F_\Ocal$ a \'et\'e
  calcul\'e par Sengupta \cite{Sengupta}.\\
  
\item[iii)] Lorsque $\Ocal$ est elliptique, $F_\Ocal$ a \'et\'e
  calcul\'e par Duflo-Vergne \cite{Duflo-Vergne1}.
\end{enumerate}

Dans \cite{pep3}, j'ai compl\'et\'e le tableau en calculant la
fonction g\'en\'eralis\'ee $F_\Ocal$ dans le cas g\'en\'eral. Je ne
vais pas donner le r\'esultat mais donner un de ses corollaires.

La fonction $\|\Phi_\Ocal\|^2$ a pour  points
critiques une $K$-orbite compos\'ee des points qui r\'ealisent le
minimum de $\|\Phi_\Ocal\|^2$: notons la $\Ccal\subset\Ocal$. Un
th\'eor\`eme (local) de forme normale d\'efinit une vari\'et\'e
symplectique $\widetilde{\Ocal}$, munie d'une action hamiltonienne de
$K$ avec application moment propre, qui est d'autre part  un fibr\'e
vectoriel au dessus de $\Ccal$. Le mod\`ele pour un voisinage de
$\Ccal$ dans $\Ocal$ est le voisinage de la section nulle du fibr\'e
vectoriel $\widetilde{\Ocal}\to\Ccal$. Notre calcul de $F_\Ocal$
montre que
$$
\mdh(\Ocal)=\mdh(\widetilde{\Ocal})\ .
$$

Ce r\'esultat est assez surprenant car le {\em mod\`ele local}
$\widetilde{\Ocal}$ se r\'ev\`ele \^etre suffisant pour calculer
l'invariant {\em global} $\mdh(\Ocal)$.

\section{Quantification g\'eom\'etrique et $K$-multiplicit\'es}
\label{sec:quant-geometrique}

Soient $(M,\omega)$ une vari\'et\'e symplectique compacte, munie de
l'action hamiltonienne d'un groupe de Lie compact connexe $K$. L'objet de la 
quantification g\'eom\'etrique est de {\em quantifier} l'action de $K$
sur $M$, c'est \`a dire de lui associer une repr\'esentation $\Qcal(M)$ du groupe $K$.
Ici nous consid\'erons des quantifications o\`u $\Qcal(M)$ est une {\em diff\'erence} 
de deux repr\'esentations de $K$: $\Qcal(M)$ appartient \`a l'anneau $R(K)$ des 
des repr\'esentations de $K$. Nous consid\'erons ici deux types de quantifications 
qui utilisent toutes deux la notion de fibr\'e en droites de Kostant-Souriau.

\medskip

 On appelle {\em fibr\'e de Kostant-Souriau} un fibr\'e en droites hermitien 
$(L,{\rm h})$ muni d'une connection hermitienne $\nabla$ dont  la courbure 
est $- i\omega$. Un tel fibr\'e existe si et seulement si la forme symplectique 
$\omega$ est {\em int\'egrale} c'est \`a dire si 
\begin{equation}
  \label{eq:condition-KS}
  \left[\frac{\omega}{2\pi}\right]\ \in\ \Hcal^2(M,\Z).
\end{equation}
Si $L$ est un fibr\'e de Kostant-Souriau, on a une action de $\kgot$ 
sur les sections de $L\to M$ donn\'ee par la formule de Kostant
\begin{equation}
  \label{eq:eq-kostant}
  \Lcal(X)=\nabla_{X_M} + i\langle\Phi,X\rangle,\quad X\in\kgot .
\end{equation}
Nous supposerons ici que l'action hamiltonienne de $K$ sur $M$ se 
rel\`eve en une action sur le fibr\'e $L$, laissant les structures 
$({\rm h},\nabla)$ invariantes.

Une des \'etapes de la quantification est le choix d'une polarisation. Ici nous 
consid\'ererons des polarisations totalement complexes, i.e. provenant d'une structure
presque complexe $J$ sur $M$, que nous ne supposerons pas int\'egrables
\`a priori.

\medskip

{\bf La quantification de type k\"ahl\'erien.}
Une structure presque complexe $K$-invariante $J$ d\'efinit un caract\`ere de 
Riemann-Roch
\begin{equation}
  \label{eq:riem-roch}
  RR^{K,J}(M,-):\K_K(M)\to R(K)\ .
\end{equation}
Ici $\K_K(M)$ d\'esigne la $\K$-th\'eorie des fibr\'es vectoriels complexes
 \'equivariants sur $M$. Si $(M,J)$ est une vari\'et\'e complexe, et 
si $E\to M$ est un fibr\'e vectoriel holomorphe, nous avons $RR^{K,J}(M,E)=
\sum_j(-1)^j \Hcal^j(M,\Ocal(E))$, o\`u $\Hcal^*(M,\Ocal(E))$ est le groupe 
de cohomologie sur $M$ du faisceau $\Ocal(E)$ des sections holomorphes de $E$.
Je donne une d\'efinition topologique du morphisme $RR^{K,J}(M,-)$ \`a la 
section \ref{subsec-loc-K}.

Pour la quantification de type k\"ahl\'erien on choisit une 
structure presque complexe $K$-invariante $J$ qui est {\em compatible} 
avec $\omega$: $\omega(-,J-)$ d\'efinit une structure riemanienne sur $M$.

\begin{defi} Une vari\'et\'e $K$-hamiltonienne compacte $(M,\omega)$ est 
$\Qcal_{\hol}$-\break 
{\em pr\'equantifi\'ee} si elle poss\`ede un fibr\'e de Kostant-Souriau $L$ 
\'equivariant. Dans ce cas, la quantification de type k\"ahl\'erien de 
l'action hamiltonienne de $K$ sur $(M,\omega)$ est  
\begin{equation}
  \label{eq:quant-hol}
  \Qcal_{\hol}(M):=RR^{K,J}(M,L)\ \in \ R(K) ,
\end{equation}
o\`u $J$ est une structure presque complexe invariante compatible avec $\omega$.
\end{defi}

\begin{rem} La quantification $\Qcal_{\hol}(M)$ ne d\'epend pas du choix du fibr\'e $L$
de Kostant-Souriau, sachant que la classe de Chern (\'equivariante) de $L$ est fix\'e.
\end{rem}

{\bf La quantification $\spinc$.} Le cadre ici a beaucoup de similitude avec ce que 
l'on appelle commun\'ement la correction m\'etaplectique. Consid\'erons une structure 
presque complexe \'equivariante $J$ {\em quelconque}, et le fibr\'e en droites complexes 
$\kappa:=\det_\C^{-1}(TM)$ correspondant. L'hypoth\`ese de travail est que 
\begin{equation}
  \label{eq:condition-spin}
  \left[\frac{\omega}{2\pi}\right]+ \frac{1}{2}{\rm c}_1(\kappa)\quad \in\ 
\Hcal^2(M,\Z) .
\end{equation}
o\`u ${\rm c}_1$ d\'esigne la premi\`ere classe de Chern. 

\begin{defi}Nous noterons $\tilde{L}$ le fibr\'e en droites complexes qui a 
pour premi\`ere classe de Chern $\left[\frac{\omega}{2\pi}\right]+ 
\frac{1}{2}{\rm c}_1(\kappa)$. Un tel fibr\'e sera appel\'e fibr\'e de 
Kostant-Souriau {\em tordu}.
\end{defi}

Si le fibr\'e $\kappa$ admet une racine carr\'ee $\kappa^{1/2}$, les conditions 
(\ref{eq:condition-KS}) et (\ref{eq:condition-spin}) co\"\i ncident et on a 
$\tilde{L}=L\otimes\kappa^{1/2}$, o\`u $L$ est un fibr\'e de Kostant-Souriau 
sur $(M,\omega)$. Dans tout les cas le produit 
\begin{equation}
  \label{eq:def-L2omega}
  L_{2\omega}:=\tilde{L}^2\otimes \kappa^{-1}
\end{equation}
est un fibr\'e de Kostant-Souriau sur $(M,2\omega)$. Comme tout \`a l'heure, on 
a une action de $\kgot$ sur les sections de $\tilde{L}\to M$. Nous parlerons de
fibr\'e de Kostant-Souriau  tordu \'equivariant lorsque l'action de $K$ sur $M$ 
se rel\`eve au fibr\'e $\tilde{L}$.

La condition (\ref{eq:condition-spin}) peut se r\'e\'ecrire en terme de structures 
$\spinc$.  La structure presque complexe $J$ 
d\'etermine une structure $\spinc$ qui a pour fibr\'e canonique le fibr\'e 
$\kappa^{-1}$. Si on tord cette structure $\spinc$ par le fibr\'e en 
droites $\tilde{L}$, on obtient une nouvelle structure $\spinc$ qui a pour 
fibr\'e canonique le fibr\'e de Kostant-Souriau $L_{2\omega}$ d\'efinie
en (\ref{eq:def-L2omega}). R\'eciproquement, 
l'existence de cette derni\`ere structure $\spinc$ implique l'existence d'un 
fibr\'e  Kostant-Souriau {\em tordu} $\tilde{L}$.
\begin{defi}\label{def-tordu} 
Une vari\'et\'e $K$-hamiltonienne compacte $(M,\omega)$ est 
$\spinc$- \break {\em pr\'equantifi\'ee} si elle poss\`ede une structure $\spinc$
\'equivariante, de fibr\'e canonique un fibr\'e de Kostant-Souriau 
sur $(M,2\omega)$. Comme on l'a bri\`evement rappel\'e cette condition est 
\'equivalente \`a l'existence d'un fibr\'e de Kostant-Souriau 
{\em tordu} \'equivariant $\tilde{L}$, associ\'e \`a une structure presque complexe
invariante $J$. Dans ce cas la quantification $\spinc$ de l'action
hamiltonienne de $K$ sur $(M,\omega)$ 
est 
\begin{equation}
  \label{eq:quant-spinc}
  \Qcal_{\spin}(M):=\esp\, RR^{K,J}(M,\tilde{L})\ \in \ R(K)  ,
\end{equation}
o\`u $\esp=\pm$ est le rapport entre les orientations de $M$ induites
par la structure presque complexe et la forme symplectique.
\end{defi}

\begin{rem}
La quantification peut \^etre aussi d\'efinie comme l'indice \'equivariant
de l'op\'erateur de Dirac associ\'e \`a la structure $\spinc$ qui a pour 
fibr\'e canonique le fibr\'e de Kostant-Souriau $L_{2\omega}$ d\'efinie
en (\ref{eq:def-L2omega}).
\end{rem}

\begin{rem}
Comme tout \`a l'heure, la quantification $\Qcal_{\spin}(M)$ 
ne d\'epend pas du choix du couple $(J,\tilde{L})$.
\end{rem}

\medskip

\subsection{La quantification commute \`a la r\'eduction}

Dans \cite{Guillemin-Sternberg82} Guillemin et Sternberg obtiennent le r\'esultat
suivant. Soient $(M,J,\omega)$ une vari\'et\'e k\"ahl\'erienne compacte telle que $\omega$ 
est enti\`ere. Soit $L$ le fibr\'e de Kostant-Souriau holomorphe associ\'e. 
Soit $K$ un groupe de Lie compact connexe agissant de mani\`ere holomorphe sur 
$L\to M$: l'action de $K$ sur $(M,\omega)$ est alors hamiltonienne, d'application 
moment $\Phi$. Supposons que $K$ agisse librement sur $\Phi^{-1}(0)$. Alors la 
vari\'et\'e r\'eduite $\Mcal_0:=\Phi^{-1}(0)/K$ est une vari\'et\'e k\"ahl\'erienne 
(lisse) munie d'une forme symplectique canonique $\omega_0$, et le fibr\'e en droites 
$\Lcal_0:=L|_{\Phi^{-1}(0)}/K$ est un fibr\'e de Kostant-Souriau holomorphe sur 
$(\Mcal_0,\omega_0)$. Guillemin et Sternberg montrent que l'op\'eration de 
restriction induit un isomorphisme
\begin{equation}
  \label{eq:iso-GS}
  \left[\Gamma_{\hol}(M,L)\right]^K\simeq \Gamma_{\hol}(\Mcal_0,\Lcal_0)\ .
\end{equation}
o\`u $[-]^K$ d\'esigne le sous-espace des vecteurs $K$-invariants. Si le 
fibr\'e $L$ est suffisamment ample, on sait d'apr\`es le th\'eor\`eme 
d'annulation de Kodaira que les autres groupes de cohomologie $\Hcal^j(M,\Ocal(L))$ 
s'annulent. L'\'equation (\ref{eq:iso-GS}) devient alors 

 \begin{equation}
  \label{eq:iso-GS-bis}
  \left[\Qcal_{\hol}(M)\right]^K= \Qcal_{\hol}(\Mcal_0)\ .
\end{equation}

Ici l'\'egalit\'e est prise dans $\Z$. Guillemin et Sternberg conjecturent alors 
que (\ref{eq:iso-GS-bis}) est encore vraie dans le cadre symplectique:
c'est alors pendant une quinzaine d'ann\'ee la conjecture intitul\'ee  
``{\em La quantification commute \`a la r\'eduction}''. 

Notons $\pi\in\f(K)$ la trace de la repr\'esentation virtuelle 
$\Qcal_{\hol}(M)$. Ce que l'on cherche ici est la valeur de l'entier 
$\left[\Qcal_{\hol}(M)\right]^K$, qui n'est autre que l'int\'egrale 
$\int_K \pi(k)dk$, o\`u $dk$ est la mesure de Haar sur $K$ normalis\'ee 
par $\int_K dk=1$.

L'id\'ee naturelle pour aborder ce probl\`eme est d'utiliser la 
formule de Atiyah-Bott-Segal-Singer \cite{Atiyah-Segal68} pour l'indice \'equivariant.
Dans la forme donn\'ee par Berline-Vergne \cite{B-G-V}, cette formule de 
l'indice donne pour tout $k\in K$
\begin{equation}
  \label{eq:ABSS}
  \pi(k e^X)=\int_{M^k}\eta_k(X),\quad X\in\{a\in\kgot,k\cdot a=a\},
\end{equation}
o\`u $\eta_k(X)$ est une 
forme caract\'eristique \'equivariante sur la sous-vari\'et\'e $M^k$ 
des points fix\'es par $k$. On remarque rapidement que (\ref{eq:ABSS})
est difficilement utilisable pour calculer l'int\'egrale $\int_K \pi(k)dk$.
Il faut en fait attendre la ``localisation non-ab\'elienne'' de Witten, pour que 
de nombreux math\'ematiciens\footnote{et math\'ematiciennes !}
s'attaquent \`a cette conjecture. Elle sera finalement 
compl\`etement d\'emontr\'ee par Meinrenken et Meinrenken-Sjamaar 
\cite{Meinrenken98,Meinrenken-Sjamaar}. Rappelons quelques faits concernant 
le travail effectu\'e pour d\'emontrer cette conjecture. On pourra trouver 
d'autres explications et r\'ef\'erences dans  
\cite{Sjamaar96,Vergne01}.

\medskip

Tout d'abord la d\'emonstration de (\ref{eq:iso-GS-bis}) s'est r\'ev\'el\'ee 
beaucoup plus difficile lorsque le groupe $K$ est non-ab\'elien. Dans le cas
ab\'elien Duistermaat, Guillemin, Meinrenken, Vergne et Wu ont 
apport\'e diff\'erentes preuves \cite{Guillemin94,Meinrenken96,DGMW,Vergne96}. 
Lorsque le groupe est non-ab\'elien, voici les travaux qui d\'emontrent
(\ref{eq:iso-GS-bis}) lorsque $0$ est une valeur r\'eguli\`ere de l'application 
moment.

\begin{itemize}
\item dans \cite{Meinrenken98} Meinrenken utilise la technique de 
coupure symplectique non-ab\'elienne d\'evelopp\'ee par Woodward,
\item dans \cite{Jeffrey-Kirwan97} Jeffrey et Kirwan d\'emontrent 
(\ref{eq:iso-GS-bis}) sous une forme asymptotique (en remplacant 
$L$ par $L^{\otimes k}$), mais en ne faisant aucune hypoth\`ese de positivit\'e
sur la structure presque complexe, 
\item dans \cite{Tian-Zhang98} Tian et Zhang d\'emontrent 
(\ref{eq:iso-GS-bis}) par une approche analytique inspir\'ee des travaux 
de Bismut et Lebeau.
\end{itemize}

Finalement Meinrenken et Sjamaar obtiennent (\ref{eq:iso-GS-bis}) dans le cas 
g\'en\'eral o\`u $0$ n'est pas forc\'ement une valeur r\'eguli\`ere de 
l'application moment \cite{Meinrenken-Sjamaar}. Dans le m\^eme temps, 
Tian et Zhang adaptent leur technique pour traiter le cas 
non-r\'egulier\footnote{N\'eanmoins ils ne traitent que le cas o\`u 
le stabilisateur g\'en\'erique de $K$ sur $M$ est fini.} \cite{Tian-Zhang-arxiv}.

Pour  donner un sens \`a (\ref{eq:iso-GS-bis}) dans le cas non-r\'egulier,
Meinrenken et Sjamaar d\'efinissent $\Qcal_{\hol}(\Mcal_0)$
au moyen d'une d\'esingularisation de $\Mcal_0$. Ils peuvent alors  
calculer les $K$-multiplicit\'es de $\Qcal_{\hol}(M)$ gr\^ace au 
``shifting trick''.

On param\`etre le dual unitaire de $K$ par l'ensemble $\Lambda^*_+$ des 
poids dominants d'une chambre de Weyl. Pour chaque $\mu\in\Lambda^*_+$ 
l'orbite coadjointe $K\cdot\mu$  munie de sa forme de Kirillov-Kostant-Souriau 
est pr\'equantifi\'ee et
$$
V_\mu:=\Qcal_{\hol}(K\cdot\mu)
$$
est la repr\'esentation irr\'eductible de $K$ de plus haut poids $\mu$.
On voit aussi que la repr\'esentation duale $V_\mu^*$ s'identifie avec
$\Qcal_{\hol}(\overline{K\cdot\mu})$, o\`u $\overline{K\cdot\mu}$
est l'orbite coadjointe munie de l'oppos\'ee de la forme symplectique 
de Kirillov-Kostant-Souriau. Pour chaque $\mu\in\Lambda^*_+$,
la multiplicit\'e de $V_\mu$ dans $\Qcal_{\hol}(M)$ est alors
\begin{eqnarray}\label{eq:shifting}
\left[\Qcal_{\hol}(M)\otimes V_\mu^*\right]^K &=&
\left[\Qcal_{\hol}(M)\otimes \Qcal_{\hol}(\overline{K\cdot\mu})\right]^K\\
&=& \left[\Qcal_{\hol}(M\times\overline{K\cdot\mu})\right]^K\nonumber\\
&=& \Qcal_{\hol}\left((M\times\overline{K\cdot\mu})_0\right)\nonumber
\end{eqnarray}
Ici la r\'eduction symplectique en $0$, $(M\times\overline{K\cdot\mu})_0$, 
s'identifie \`a la r\'eduction symplectique 
$\Mcal_\mu:=\Phi^{-1}(\mu)/K_\mu$ qui est $\Qcal_{\hol}$-pr\'equantifi\'ee par 
$\Lcal_\mu:=\left(L|_{\Phi^{-1}(\mu)}\otimes \C_{-\mu}\right)/K_\mu$.  
Le th\'eor\`eme ``La quantification commute \`a la r\'eduction'' prend 
alors la forme finale suivante.

\begin{theo}[\cite{Meinrenken98,Meinrenken-Sjamaar}] \label{th:Q-R-hol}
Pour toute vari\'et\'e $K$-hamiltonienne $(M,\omega)$ \break 
$\Qcal_{\hol}$-pr\'equantifi\'ee, on a
\begin{equation}
\Qcal_{\hol}(M)=\sum_{\mu\in \Lambda^*_+}
\Qcal_{\hol}(\Mcal_\mu)\, V_\mu\ .
\end{equation}
\end{theo}

\bigskip

Mon travail dans le domaine de la quantification a commenc\'e en 1999 (un peu 
apr\`es la bataille). Dans \cite{Vergne96}, Vergne reprend une id\'ee d'Atiyah de 
localisation  $\K$-th\'eorique, pour d\'emontrer (\ref{eq:iso-GS-bis}) dans 
le cas d'une action du cercle. J'ai d\'evelopp\'e ce proc\'ed\'e de 
localisation pour les actions 
d'un groupe de Lie compact connexe quelconque \cite{pep4}. Cette technique me permet
alors de red\'emontrer le th\'eor\`eme \ref{th:Q-R-hol}, mais aussi
d'apporter une g\'en\'eralisation dans un cadre non-symplectique.

Dans \cite{pep5}, je reprends cette m\'ethode de localisation 
$\K$-th\'eorique pour d\'emontrer que le slogan 
``la quantification commute \`a la r\'eduction'' est vrai 
dans le cadre de la quantification $\spinc$. Dans ce m\^eme article je montre 
que ce r\'esultat est encore valable dans certains cas o\`u la 
vari\'et\'e $M$ est non-compacte: ce sont les orbites coadjointes qui param\`etrent
la {\em s\'erie discr\`ete} d'un groupe de Lie r\'eel semi-simple.

\medskip

Dans la prochaine sous-section, je r\'esume les diff\'erents r\'esultats obtenus.
Ensuite, je donne un bref apercu sur ce proc\'ed\'e de localisation 
$\K$-th\'eorique.

\subsection{R\'esultats sur la  quantification de type k\"ahl\'erien.}
J'ai red\'emontr\'e dans \cite{pep5}, le th\'eor\`eme 
``la quantification commute  \`a la r\'eduction'' dans le cadre 
hamiltonien, lorsque $0$ n'est pas forc\'ement une valeur r\'eguli\`ere
de l'application moment. N\'eanmoins l'essentiel de cet article est 
consacr\'e \`a l'\'etude de la validit\'e de (\ref{eq:iso-GS-bis}) en dehors 
du cadre symplectique. Soit $M$ une vari\'et\'e diff\'erentiable 
munie 
\begin{itemize}
\item d'une action d'un groupe de Lie compact connexe $K$,
\item d'une structure presque complexe $J$ qui est $K$-invariante, 
\item d'une $2$-forme ferm\'ee, $K$-invariante et {\em int\'egrable}. 
\end{itemize}
Contrairement \`a la section pr\'ec\'edente, on ne fait {\em aucune hypoth\`ese} sur
le caract\`ere {\em non-d\'eg\'en\'er\'e} de la $2$-forme en question, ni sur une 
{\em compatibilit\'e} de $J$ avec la $2$-forme. Soit $L\to M$ le fibr\'e 
en droites complexes associ\'e \`a cette $2$-forme: cette derni\`ere est not\'ee 
dor\'enavant $\omega^L$. {\em On suppose dans cette partie que l'action de 
$K$ sur $M$ se rel\`eve au fibr\'e vectoriel $L$.}

On munit le fibr\'e $L$ d'une structure hermitienne  et d'une connexion 
hermitienne $\nabla^L$, tous deux $K$-invariants, tels que $(\nabla^L)^2=-i\omega^L$. 
Ces donn\'ees d\'eterminent une application
\'equivariante $\Phi^L:M\to \kgot^*$ satisfaisant
\begin{equation}\label{eq.moment.L}
\Lcal^L(X)-\nabla^L_{X_M}=\imath \langle \Phi^L,X\rangle,\quad 
X\in\kgot\ .
\end{equation}
La formule de Bianchi \'equivariante (voir Prop. 7.4 dans \cite{B-G-V})
donne 
$$
d\langle \Phi^{L},X\rangle=-\omega^{L}(X_{M},-).
$$
Ainsi l'application $\Phi^L$ est une {\em application moment abstraite} 
au sens de Karshon: elle est \'equivariante, et pour tout $X\in\kgot$, la 
fonction $\langle \Phi^L,X\rangle$ is locallement constante sur la sous-vari\'et\'e 
$M^X$ \cite{Karshon.98}. On remarque ici que la donn\'ee de $(L,\omega^L,\Phi^L)$ 
est une g\'en\'eralisation de la quantification de Kostant-Souriau au cadre 
non-symplectique. Comme tout \`a l'heure, on s'int\'eresse \`a la
multiplicit\'e $\left[RR^{K,J}(M,L)\right]^K\in\Z$. 

\medskip

Nous supposons ici que $0$ est {\em une valeur r\'eguli\`ere de $\Phi^L$}: 
alors $\Zcal:=(\Phi^L)^{-1}(0)$ est une sous-vari\'et\'e de $M$ munie d'une action 
(localement) libre de $K$. On consid\`ere le quotient 
$\Mcal_0:=\Zcal/K$ qui est une $V$-vari\'et\'e diff\'erentiable munie du 
$V$-fibr\'e en droites complexes $L_0:=(L|_\Zcal)/K$.

\begin{theo}[\cite{pep5}]\label{theo:Q-R-1} 
La structure presque complexe $J$ d\'etermine 
une structure $\spinc$ sur $\Mcal_0$: soit  $\Qcal(\Mcal_0,-): 
\K(\Mcal_0)\to \Z$ le morphisme d\'etermin\'e par cette 
structure $\spinc$. On a
\begin{equation}\label{eq.theorem.A}
 \left[ RR^{K,J}(M,L^{\otimes k})\right]^K=
 \Qcal\Big(\Mcal_0,(L_0)^{\otimes k}\Big), \quad 
 k\in\N-\{0\},
\end{equation}
si l'une des deux conditions est satisfaite:

(i) $K=T$ is a tore,

(ii) $k\in \N$ est assez grand, de telle mani\`ere que la boule
$\{ \xi\in\kgot^{*},\, \parallel\xi\parallel\leq\frac{1}{k}
\parallel\theta\parallel \}$ est contenue dans l'ensemble des valeurs 
r\'eguli\`eres de $\Phi^L$. Ici $\theta$ est la somme des racines positives
de $K$ et $\parallel\cdot\parallel$ est une norme euclidienne $K$-invariante 
sur $\kgot^*$.
\end{theo}

Ce th\'eor\`eme montre que le slogan ``la quantification commute
\`a la r\'eduction'' est un ph\'enom\`ene assez g\'en\'eral lorsque 
$0$ est une valeur r\'eguli\`ere: il est toujours vrai dans le cas ab\'elien, 
et dans le cas g\'en\'eral il est vrai de fa\c con  asymptotique. Dans le 
cas ab\'elien, un r\'esultat similaire\footnote{Leurs hypoth\`eses de travail 
sont n\'eanmoins diff\'erentes: ils supposent que le lieu des points fixes est 
fini, et travaille avec une  structure presque complexe stable.}
a \'et\'e obtenu par Ginzburg, Guillemin et Karshon 
\cite{GGK-bouquin}. 

\medskip

L'\'equation (\ref{eq.theorem.A}) peut \^etre r\'e\'ecrite lorsque 
$J$ induit une  structure presque complexe $J_0$ on $\Mcal_0$. 
C'est le cas lorsque on a la d\'ecomposition
\begin{equation}\label{eq.J.red}
\T M\vert_{\Zcal}=\T\Zcal\oplus J(\kgot_{\Zcal})\quad {\rm avec}\quad 
\kgot_{\Zcal}:=\{X_{\Zcal},\, X\in\kgot\}\ .
\end{equation}
Cette derni\`ere condition est toujours satisfaite dans le cadre hamiltonien 
lorsque $J$ est compatible avec la forme symplectique. La 
condition (\ref{eq.J.red}) appara\^\i t d\'ej\`a dans les travaux de 
Jeffrey-Kirwan \cite{Jeffrey-Kirwan97}, and Cannas da Silva-Karshon-Tolman 
\cite{C-K-T}. Consid\'erons pour tout $z\in\Zcal$, l'endomorphisme
$\Dcal_z:\kgot\to\kgot^*$ d\'efinit par $\Dcal_z(X)=-
d\Phi^L|_z(J(X_{\Zcal}|_z))$. Alors (\ref{eq.J.red}) est \'equivalent 
\`a $\det(\Dcal_z)\neq 0,\ \forall z\in\Zcal$. Soit $J_{\Dcal}$ la structure 
complexe sur le fibr\'e trivial  $(\kgot\oplus\kgot^*)\times\Zcal$ 
d\'efinie par la matrice
$J_{\Dcal}:=\left(\begin{array}{cc}0   & -\Dcal^{-1}  \\ \Dcal & 0
\end{array}\right)$. On une autre struture complexe $J_\C$ sur 
$(\kgot\oplus\kgot^*)\times\Zcal$, d\'efinie par $J_{\C}:=
\left(\begin{array}{cc}0   & -I^{-1}  \\ I & 0\end{array}\right)$, 
o\`u $I:\kgot\stackrel{\sim}{\to}\kgot^*$ est 
un isomorphisme induit par un produit scalaire invariant sur $\kgot$. 
Le fibr\'e trivial $(\kgot\oplus\kgot^*)\times\Zcal$ est munit de deux 
modules de cliffords irr\'eductibles, 
$$
\wedge^{\bullet}_{J_\Dcal}(\kgot\oplus\kgot^*)\times\Zcal\quad
{\rm et}\quad
\wedge^{\bullet}_{J_\C}(\kgot\oplus\kgot^*)\times\Zcal \ .
$$
On note $\Lfibre\to \Mcal_0$ le fibr\'e en droites complexes  
d\'efini par
$$
\pi^*(\Lfibre):={\rm Hom}_{_{Clifford}}
\Big(
\wedge^{\bullet}_{J_\C}(\kgot\oplus\kgot^*)\times\Zcal \, , \, 
\wedge^{\bullet}_{J_\Dcal}(\kgot\oplus\kgot^*)\times\Zcal
\Big)\ .
$$
Ici $\pi:\Zcal\to\Mcal_0$ d\'esigne le quotient par $K$.

\begin{prop}[\cite{pep5}]
Si l'hypoth\`ese (\ref{eq.J.red}) est v\'erifi\'ee, 
l'\'egalit\'e (\ref{eq.theorem.A}) devient
\begin{equation}\label{eq.theorem.A.bis}
 \left[ RR^{K,J}(M,L^{\otimes k})\right]^K= \pm
 RR^{J_0}\Big(\Mcal_0,(L_0)^{\otimes k}\otimes\Lfibre\Big).
\end{equation}
o\`u  $\pm$ est le signe du d\'eterminant de $\Dcal$.
\end{prop}

Le th\'eor\`eme \ref{theo:Q-R-1}, sous la forme donn\'ee par 
(\ref{eq.theorem.A.bis}), \'etend le r\'esultat principal de Jeffrey-Kirwan
dans \cite{Jeffrey-Kirwan97} au cadre non-symplectique. Elles obtiennent l'\'egalit\'e
(\ref{eq.theorem.A.bis}) dans le cadre hamiltonien lorsque $J$ n'est pas suppos\'e
 compatible avec la forme symplectique (il semble n\'eanmoins 
qu'elles n\'egligent le facteur $\Lfibre$).

\bigskip

{\em Fibr\'es $\Phi$-positifs}

\medskip

On peut travailler en sens inverse que pr\'ec\'edemment. Supposons que 
$\Phi:M\to \kgot^*$ soit une application moment abstraite au sens de 
Karshon: $\Phi$ est \'equivariante, et pour tout $X\in\kgot$, la 
fonction $\langle \Phi,X\rangle$ est locallement constante sur la 
sous-vari\'et\'e $M^X$. Supposons de plus que $0$ est une 
valeur r\'eguli\`ere de $\Phi$. Comme tout \`a l'heure on a une 
fibr\'e principal $\Zcal\to\Mcal_0:=\Zcal/K$ avec $\Zcal=\Phi^{-1}(0)$,
la structure presque complexe $J$ induit une structure $\spinc$ sur
$\Mcal_0$ et on a donc un morphisme $\Qcal(\Mcal_0,-):\K(\Mcal_0)\to \Z$. 
Pour quels fibr\'es \'equivariants la relation (\ref{eq.theorem.A}) 
est-elle encore vraie ? 

Soit $E$ un fibr\'{e} vectoriel  hermitien $K$-\'{e}quivariant sur $M$.  
Si $X\in \kgot$ et si $m\in M$ est un z\'{e}ro du champ $X_M$, alors le groupe 
\`{a} un param\`{e}tre $\exp tX$ op\`{e}re dans la fibre  $E_m$ du fibr\'{e} 
$E$ au point $m$. On note $\Lcal^E_m(X)$ l'application infinit\'{e}simale 
de $X$ dans $E_m$. Alors $i\Lcal^E_m(X)$ est un endomorphisme 
hermitien de $E_m$. 

\begin{defi}\label{def:stric-positif}
Soit $E\to M$ un fibr\'e hermitien $K$-\'equivariant.
\begin{itemize}
  \item  On dit que $E$ est  $\Phi$-positif si pour tout $X\in \kgot$ et 
tout z\'{e}ro  $m$ du champ  $X_M$,  l' op\'{e}rateur 
$-i\langle\Phi(m),X\rangle\Lcal^E_m(X)$ a  toutes ses valeurs 
propres  positives ou nulles.   
  \item  On dit que $E$ est strictement  $\Phi$-positif,  si pour 
tout $X\in \kgot$ et tout z\'{e}ro  $m$ du champ  $X_M$, tel que 
$\langle\Phi(m), X\rangle\neq 0$, l'op\'{e}rateur  
$-i\langle\Phi(m),X\rangle\Lcal^E_m(X)$  a  toutes ses valeurs propres  
strictement positives. 
\end{itemize}
\end{defi}

Tian-Zhang introduisent une notion similaire dans le cadre hamiltonien 
\cite{Tian-Zhang98}. Tout fibr\'e en droites \'equivariant $L$ est strictement 
positif par rapport \`a l'application moment  abstraite $\Phi^L$  
d\'etermin\'ee par (\ref{eq.moment.L}). Le fibr\'e trivial $M\times\C\to M$ 
avec une action triviale de $K$ sur $\C$ est un fibr\'e $\Phi$-positif.  
 
\begin{theo}[\cite{pep5}] Soit $J$ une structure presque
complexe \'equivariante sur $M$. 
Pour tout fibr\'e $E$ qui est $\Phi$-{\em strictement positif}, on a
\begin{equation}\label{eq.theorem.B}
 \left[ RR^{K,J}(M,E^{\otimes k})\right]^K=
 \Qcal\Big(\Mcal_0,(E_0)^{\otimes k}\Big), \quad 
 k\in\N-\{0\},
\end{equation}
si l'une des deux conditions est satisfaite:

(i) $G=T$ is a tore,

(ii) $k\in \N$ est assez grand.

L'\'egalit\'e (\ref{eq.theorem.B}) est encore vraie pour 
les fibr\'es {\em positifs}, lorsque on est dans la cadre hamiltonien:
$(M,\omega)$ est une vari\'et\'e symplectique, $\Phi$ est l'application moment
associ\'ee \`a l'action hamiltonienne de $K$ sur $(M,\omega)$,
 $J$ est compatible avec $\omega$, et $0\in\Phi(M)$.
\end{theo}

Ce th\'eor\`eme g\'en\'eralise un r\'esultat de Tian-Zhang obtenu dans le cadre
hamiltonien \cite{Tian-Zhang98}.

\medskip

Lorsque $0$ est une valeur r\'eguli\`ere 
``la quantification commute  \`a la r\'eduction'' est 
un  ph\'enom\`ene assez g\'en\'eral. La diff\'erence entre
le cadre g\'en\'eral ``d'application moment abstraite'' et la cadre hamiltonien
est  nette lorsque $0$ n'est plus une valeur r\'eguli\`ere. Meinrenken
et Sjamaar ont les premiers montr\'e que dans le cadre hamiltonien 
l'\'egalit\'e (\ref{eq:iso-GS-bis}) est encore vraie quitte \`a remplacer $0$ 
par une valeur r\'eguli\`ere proche de $0$. Ceci ne marche plus dans le cadre 
non-symplectique (voir \cite{GGK-bouquin} pour une discussion sur ce probl\`eme 
et quelques \'el\'ements de r\'eponse).

\subsection{R\'esultats sur la  quantification $\spinc$}

Avant mon papier \cite{pep5}, le cadre de la quantification $\spinc$ a d\'ej\`a \'et\'e
abord\'e dans les travaux de Cannas da Silva-Karshon-Tolman \cite{C-K-T} 
et de Vergne \cite{Vergne96} lorsque le
groupe est ab\'elien, et ceux de Jeffrey-Kirwan \cite{Jeffrey-Kirwan97} 
pour un groupe compact connexe g\'en\'eral. 

Les r\'esultats de Jeffrey-Kirwan expriment
la multiplicit\'e $[\Qcal_{\spin}(M)]^K$ au moyen de la vari\'et\'e r\'eduite
$\Mcal_0=\Phi^{-1}(0)/K$. Elle obtiennent ce r\'esultat sous la premi\`ere 
hypoth\`ese qu'une boule centr\'ee en $0$ de taille assez grande est incluse 
dans l'ouvert des valeurs r\'eguli\`eres de $\Phi$, et sous la seconde 
hypoth\`ese que la d\'ecomposition (\ref{eq.J.red}) existe. Ces r\'esultats 
souffrent du fait que Jeffrey et Kirwan traitent la situation $\spinc$ de mani\`ere
analogue \`a la quantification de type k\"ahl\'erien. Montrons pourquoi on doit 
en fait proc\'eder diff\'eremment. 

\medskip

{\em Quels mod\`eles ?} Revenons au proc\'ed\'e du ``shifting  trick'' que 
l'on a d\'ej\`a expliqu\'e (voir (\ref{eq:shifting})). Consid\'erons tout 
d'abord la quantification de type k\"ahl\'erien d'une vari\'et\'e hamiltonienne 
quantifi\'ee $(M,\omega,\Phi)$. Pour calculer les $K$-multiplicit\'es de
$\Qcal_{\hol}(M)$, on a recours \`a des {\em mod\`eles}. Pour chaque poid
dominant $\mu\in\Lambda^*_+$, l'orbite coadjointe $K\cdot\mu$ est l'unique 
orbite coadjointe $\Ocal\subset\kgot^*$ $\Qcal_{\hol}$-pr\'equantifi\'ee telle que  
$\Qcal_{\hol}(\Ocal)$ soit \'egal \`a la repr\'esentation irr\'eductible 
de plus haut poid $\mu$ (que l'on note $V_\mu$). Pour la quantification de 
type k\"ahl\'erien, l'orbite coadjointe $K\cdot\mu$ est ainsi {\em le mod\`ele} 
de la repr\'esentation irr\'eductible $V_\mu$.

\medskip

Quels sont les mod\`eles pour la quantification $\spinc$ ? Notons $\rho_c$ la
demi-somme des racines positives. Consid\'erons pour chaque poid
dominant $\mu\in\Lambda^*_+$, l'orbite coadjointe $K\cdot(\mu+\rho_c)$ qui est 
une orbite r\'eguli\`ere, i.e. isomorphe \`a $K/T$. Si on munit 
$K\cdot(\mu+\rho_c)$ de sa forme symplectique de Kirillov-Kostant-Souriau et
de l'unique structure complexe invariante compatible,
on voit ais\'ement que $K\times_T\C_\mu\to K/T\simeq K\cdot(\mu+\rho_c)$ est 
un fibr\'e de Kostant-Souriau {\em tordu} (voir d\'efinition
\ref{def-tordu}). Et dans ce cas on a
$$
\Qcal_{\spin}\Big(K\cdot(\mu+\rho_c)\Big)=V_\mu \ .
$$
L'orbite coadjointe $K\cdot(\mu+\rho_c)$ est, pour la 
quantification $\spinc$, le mod\`ele de la repr\'esentation 
irr\'eductible $V_\mu$. Nous montrons dans \cite{pep5} que 
le slogan ``la quantification commute \`a la r\'eduction'' est vrai
dans le cadre $\spinc$.

\begin{theo}[\cite{pep5}]\label{th:Q-R-spinc} 
Soit $(M,\omega,\Phi)$ une vari\'et\'e $K$-hamiltonienne 
$\spinc$ \break pr\'equantifi\'ee par un fibr\'e de Kostant-Souriau
tordu \'equivariant $\tilde{L}$. On note $L_{2\omega}$ le fibr\'e de Kostant-Souriau
sur $(M,2\omega)$ correspondant (voir (\ref{eq:def-L2omega})).

$\bullet$ {\bf La pr\'equantification $\spinc$ est pr\'eserv\'ee par les 
r\'eductions symplectiques}. Soit $\mu\in\Lambda^*_+$ un poids dominant
tel que $\mu+\rho_c$ est une valeur r\'eguli\`ere de $\Phi$. La 
$V$-vari\'et\'e r\'eduite $\Mcal_{\mu+\rho_c}=\phi^{-1}(\mu+\rho_c)/T$ est alors 
$\spinc$-pr\'equantifi\'ee par le fibr\'e de Kostant-Souriau
tordu $\left(\tilde{L}|_{\phi^{-1}(\mu+\rho_c)}\times\C_{-\mu}\right)/T$: on note
$\Qcal_{\spin}(\Mcal_{\mu+\rho_c})\in\Z$ la quantification correspondante.
Dans le cas g\'en\'eral o\`u $\mu+\rho_c$ n'est pas forc\'ement 
une valeur r\'eguli\`ere de $\Phi$, on choisit une valeur 
r\'eguli\`ere $\xi\in\tgot^*$ suffisamment proche de $\mu+\rho_c$.
La vari\'et\'e r\'eduite $\Mcal_{\xi}$ est munie alors d'une structure
$\spinc$ qui a pour fibr\'e canonique 
$\left(L_{2\omega}|_{\phi^{-1}(\xi)}\times\C_{-\mu}\right)/T$. L'indice de 
l'op\'erateur de Dirac-$\spinc$ correspondant ne d\'epend pas de $\xi$:
on le note encore $\Qcal_{\spin}(\Mcal_{\mu+\rho_c})\in\Z$.

\medskip

$\bullet$ {\bf la quantification commute \`a la r\'eduction}. 
Si les stabilisateurs pour l'action de $K$ sur $M$ sont {\bf ab\'eliens}
on a
\begin{equation}
  \label{eq-QR-spinc}
\Qcal_{\spin}(M)=\sum_{\mu\in\Lambda^*_+}
\Qcal_{\spin}(\Mcal_{\mu+\rho_c})\, V_\mu\quad {\rm dans}\quad R(K) .
\end{equation}
\end{theo}

\medskip

La condition ``stabilisateurs ab\'eliens'' est ici n\'ecessaire. Cela vient du fait
qu'il n'y pas unicit\'e des mod\`eles. Consid\'erons par exemple la 
repr\'esentation triviale $V_0$. Comme on l'a vu l'orbite coadjointe 
$K\cdot\rho_c$ est $\spinc$-pr\'equantifi\'ee et satisfait
$\Qcal_{\spin}(K\cdot\rho_c)=V_0$. Mais d'autre orbites coadjointes satisfont
les m\^emes conditions. On a l'orbite  $\{0\}$, mais aussi 
toutes les orbites $K\cdot(\rho_c-\rho_{c,\sigma})$ o\`u $\sigma$ est une face 
ouverte de la chambre de Weyl et $\rho_{c,\sigma}$ est la demi-somme des 
racines positives qui s'annulent sur $\sigma$. Ici le fibr\'e trivial
$K\cdot(\rho_c-\rho_{c,\sigma})\times \C \to K\cdot(\rho_c-\rho_{c,\sigma})$
est un fibr\'e de Kostant-Souriau tordu et  
$\Qcal_{\spin}\left(K\cdot(\rho_c-\rho_{c,\sigma})\right)=V_0$.

Par contre l'orbite $K\cdot\rho_c$ est le seul mod\`ele qui satisfait
la condition ``stabilisateurs ab\'eliens''. Ceci explique la nature de cette
condition.

\medskip

Dans la prochaine section, je rappelle le r\'esultat principal de \cite{pep5} qui
peut \^etre r\'esum\'e ainsi.  Les orbites
coadjointes d'un groupe r\'eel semi-simple $G$ qui param\`etrent la s\'erie 
discr\`ete de $G$ sont des exemples {\bf non-compacts} de quantification
$\spinc$ o\`u le slogan ``la quantification commute \`a la r\'eduction'' 
est vrai.

\subsection{La s\'erie discr\`ete}\label{subsec:serie}

Soit $G$ un groupe r\'eel semi-simple connexe de centre fini. Par d\'efinition,
la s\'erie discr\`ete de $G$ est l'ensemble (modulo isomorphisme) 
des repr\'esentations unitaires irr\'eductibles de $G$ qui sont de carr\'e 
int\'egrable. Soit $K$ un groupe compact maximal de $G$, et $T$ un tore maximal
de $K$. Harish-Chandra a montr\'e que la s\'erie discr\`ete de $G$ est non-vide
si et seulement si $T$ est un sous-groupe de cartan de $G$. Dans cette section,
on suppose donc que $T$ est un sous-groupe de Cartan de $G$. Soient 
$\ggot,\kgot,\tgot$ les alg\`ebres de Lie de $G,K,T$, et $\ggot^*,\kgot^*,
\tgot^*$ leur dual. Un \'el\'ement $X\in\tgot$ est dit $G$-{\em r\'egulier} si
le sous groupe stabilisateur $G_X$ est \'egal \`a $T$. 
Soit $\Lambda^*\subset\tgot^*$ l'ensemble des poids (r\'eels) de $T$.
On note respectivement $\Rgot_c\subset\Rgot\subset\Lambda$ l'ensemble des 
racines pour l'action de $T$ sur $\kgot\otimes\C$ et $\ggot\otimes\C$. 
On fixe une chambre de Weyl $\tgot^*_+$ pour le couple $(K,T)$. 

Harish-Chandra param\`etre \cite{Harish-Chandra} la s\'erie discr\`ete de $G$ par un sous
ensemble discret $\widehat{G}_d$ compos\'e d'\'el\'ements $G$-r\'eguliers de $\tgot^*_+$:
pour $\lambda\in\widehat{G}_d$, notons $\Hcal_\lambda$ la repr\'esentation 
de $G$ correspondante.
On a en fait un moyen g\'eom\'etrique de d\'efinir $\widehat{G}_d$:
``un \'el\'ement $G$-r\'egulier $\lambda$ appartient \`a $\widehat{G}_d$ si et 
seulement si l'orbite coadjointe $G\cdot\lambda$ est $\spinc$-pr\'equantifi\'ee''.

Fixons pour le reste de cette section $\lambda\in\widehat{G}_d$ et
rappelons la r\'ealisation donn\'ee par Schmid \cite{Schmid} 
de la repr\'esentation $\Hcal_\lambda$ comme $\spinc$-quantification de 
l'orbite coadjointe $G\cdot\lambda$. 

Soit $\Rgot^{+,\lambda}\subset\Rgot$ 
l'ensemble de racines positives d\'efini par $\lambda$: 
$\alpha\in\Rgot^{+,\lambda}\Longleftrightarrow (\alpha,\lambda)>0$.
Soit $\rho$ la demi-somme des \'el\'ements de $\Rgot^{+,\lambda}$.
On travaille avec la structure complexe $J$ sur $G\cdot\lambda$,  qui est 
$G$-invariante, et qui est d\'etermin\'ee par la condition: 
$\alpha\in\Lambda^*$ est un poids de l'action de $T$ sur l'espace tangent
$(\T_\lambda(G\cdot\lambda),J)$ si et seulement si $\alpha\in\Rgot^{+,\lambda}$.
La condition $\lambda\in\widehat{G}_d$ impose que $\lambda-\rho$ est un poids 
de $T$, et on voit alors que le fibr\'e en droites
$$
\tilde{L}:=G\times_T\C_{\lambda-\rho}
$$
est un fibr\'e de Kostant-Souriau tordu sur $G\cdot\lambda\simeq G/T$.
Rappelons que $G\cdot\lambda$ est munie de la forme symplectique de 
Kirillov-Kostant-Souriau, l'action de $G$ sur $G\cdot\lambda$ est hamiltonienne
avec pour application moment l'inclusion $G\cdot\lambda\croc\ggot^*$.

Schmid r\'ealise $\Hcal_\lambda$ comme un espace de cohomologie $\lde$ \`a valeurs
dans $\tilde{L}$ \cite{Schmid}. Le fibr\'e en droites complexes $\tilde{L}$ poss\`ede une 
structure holomorphe canonique. Soit
$\Omega^k(\tilde{L})$ l'espace des formes diff\'erentielles sur $G\cdot\lambda$, 
de type $(0,k)$, et \`a valeurs dans $\tilde{L}$. Soit
$\overline{\partial}_{\tilde{L}}:\Omega^k(\tilde{L})\to 
\Omega^{k+1}(\tilde{L})$ l'op\'erateur de Dolbeault. Le choix de m\'etriques 
$G$-invariantes sur $G\cdot\lambda$ et $\tilde{L}$ permet de d\'efinir
l'op\'erateur $\overline{\partial}^*_{\tilde{L}}$ (adjoint formel de 
$\overline{\partial}_{\tilde{L}}$) et l'op\'erateur de Dolbeault-Dirac 
$ \overline{\partial}_{\tilde{L}}+ \overline{\partial}^*_{\tilde{L}}$. 

La cohomologie $\lde$ de $\tilde{L}$, que l'on note 
$\ch^*_{(2)}(G\cdot\lambda,\tilde{L})$,  est par d\'efinition le noyau
de l'op\'erateur $\overline{\partial}_{\tilde{L}}+ 
\overline{\partial}^*_{\tilde{L}}$ sur le sous-espace de
$\Omega^*(\tilde{L})$ des \'el\'ements de carr\'e int\'egrable.

\medskip

\begin{theo}(Schmid).\ Soit $\lambda\in\widehat{G}_{d}$.

(i) Si  $k\neq\frac{\dim(G/K)}{2}$, alors $\ \ch^k_{(2)}(G\cdot\lambda,\tilde{L})=0$ .

(ii) Si $k=\frac{\dim(G/K)}{2}$, alors
$\ \ch^k_{(2)}(G\cdot\lambda,\tilde{L})=\Hcal_{\lambda}$.
\end{theo}

\medskip

Ainsi la repr\'esentation $\Hcal_{\lambda}$ est la quantification 
$\spinc$ de l'action de $G$ sur l'orbite coadjointe $G\cdot\lambda$ 
comme indice $\lde$ de l'op\'erateur de Dolbeault-Dirac 
$\overline{\partial}_{\tilde{L}}+ \overline{\partial}^*_{\tilde{L}}$: ici 
$(-1)^{\frac{\dim(G/K)}{2}}$ est le rapport entre les orientations de 
$G\cdot\lambda$ induites par la structure complexe $J$ et la forme 
symplectique. La restriction $\Hcal_\lambda|_K$ est 
la quantification $\spinc$ de l'action de $K$ sur l'orbite coadjointe 
$G\cdot\lambda$. La repr\'esentation $\Hcal_\lambda|_K$ admet une d\'ecomposition 
en repr\'esentation irr\'eductibles de $K$
$$
\Hcal_{\lambda}\vert_{K}=\sum_{\mu\in\Lambda^*_+}\mm_{\mu}(\lambda)\, 
V_{\mu}\ ,
$$
o\`u les multiplicit\'es $\mm_{\mu}(\lambda)$ satisfont les formules combinatoires 
dites ``formules de Blattner'' \cite{Hecht-Schmid}. L'objet de l'article \cite{pep5}
est de montrer que ces $K$-multiplicit\'es satisfont le principe
de Guillemin-Sternberg. L'action de $K$ sur $G\cdot\lambda$ est hamiltonienne avec 
pour application moment $\Phi:G\cdot\lambda\to\kgot$ le compos\'e de
l'inclusion $G\cdot\lambda\croc\ggot^*$ avec la projection $\ggot^*\to\kgot^*$.
L'orbite coadjointe $G\cdot\lambda$ est non-compacte, mais l'application 
$\Phi$ est {\em propre} \cite{pep3}. Ainsi les vari\'et\'es r\'eduites 
$\Phi^{-1}(\xi)/K_\xi$ 
sont compactes. La premi\`ere partie du th\'eor\`eme \ref{th:Q-R-spinc} 
s'applique: la pr\'equantification $\spinc$ de 
$G\cdot\lambda$ induit pour tout $\mu\in\Lambda^*_+$ une pr\'equantification 
$\spinc$ sur chaque r\'eduction symplectique $(G\cdot\lambda)_{\mu+\rho_c}:=
\Phi^{-1}(\mu+\rho_c)/T$. Comme ces vari\'et\'es r\'eduites sont compactes, 
on peut d\'efinir\footnote{Par d\'esingularisation 
si n\'ecessaire.} la quantit\'e $\Qcal_{\spin}((G\cdot\lambda)_{\mu+\rho_c})\in\Z$ comme
l'indice de l'op\'erateur $\spinc$-Dirac associ\'e. Je montre alors que la 
deuxi\`eme partie du th\'eor\`eme \ref{th:Q-R-spinc} est encore vraie dans ce cas.
Pour cela j'utilise les formules combinatoires de Blattner.

\begin{theo}[\cite{pep5}]\label{theo-serie-discrete}
\begin{equation}
  \label{eq-Glambda-spinc}
\Hcal_\lambda|_K=\sum_{\mu\in\Lambda^*_+}
\Qcal_{\spin}((G\cdot\lambda)_{\mu+\rho_c})\, V_\mu \  .
\end{equation}
\end{theo}

\subsection{Localisation $\K$-th\'eorique}\label{subsec-loc-K} 
Soit $M$ une vari\'et\'e diff\'erentiable compacte munie d'une action d'un groupe
de Lie compact connexe $K$. Supposons $M$ munie d'une structure presque complexe 
$K$-invariante $J$. Soit $\K_K(M)$ la 
$\K$-th\'eorie des fibr\'es vectoriels complexes \'equivariants sur $M$.
La structure presque complexe permet de d\'efinir  le caract\`ere de Riemann-Roch
\begin{equation}
  \label{eq:RR}
  RR^{K,J}(M,-): \K_K(M) \longrightarrow  R(K)\\
\end{equation}

\medskip

Donnons la d\'efinition topologique de ce morphisme. 
Soit $p:\T M\to M$ la projection, et $E^\pm\to M$ des fibr\'es vectoriels complexes
\'equivariants. Une section $K$-\'equivariante 
$\sigma \in \Gamma(\T M,\hom(p^{*}E^{+},p^{*}E^{-}))$ est appel\'e un {\em symbole.}
L'ensemble des $(m,v)\in \T M$ o\`u  
$\sigma(m,v): E^{+}_{m}\to E^{-}_{m}$ n'est pas inversible est 
{\em l'ensemble caract\'eristique} of $\sigma$, que l'on note $\Char(\sigma)$. 
Un symbole $\sigma$ est dit elliptique si $\Char(\sigma)$ est compact: il d\'efinit
alors un \'el\'ement de la $\K$-th\'eorie \'equivariante de $\T M$ \`a support 
compact, qui est not\'ee $\K_K(\T M)$. Nous avons dans ce cas un morphisme
$\indice^K_M:\K_K(\T M)\to R(K)$ 
\cite{Atiyah-Segal68,Atiyah-Singer-1}.

La structure presque complexe d\'efinit le symbole 
$$
\Thom_K(\T M)\in\Gamma\left(M,\hom(p^{*}(\wedge_{\C}^{pair} \T M),\,p^{*}
(\wedge_{\C}^{impair} \T M))\right)\ .
$$
En $(m,v)\in \T M$, le morphisme $\Thom_K(\T M)(m,v):\wedge_{\C}^{pair} \T_m M
\longrightarrow \wedge_{\C}^{impair} \T_m M$ correspond \`a l'action de Clifford de $v$.
L'ensemble caract\'eristique de $\Thom_K(\T M)$ est $M$ qui est compact: le symbole de
$\Thom_K(M)$ est elliptique. Pour tout $K$-fibr\'e vectoriel complexe $E$, on peut 
consid\'erer le ``produit'' $\Thom_K(\T M)\otimes p^*(E)$ qui est encore un symbole
elliptique. L'isomorphisme de Bott-Thom affirme que l'application
\begin{eqnarray}
  \label{eq:iso-thom}
   \Thom_J:\K_K(M)& \longrightarrow & \K_K(\T M)\\
     E & \longmapsto & \Thom_K(\T M)\otimes p^*(E)\ , \nonumber
  \end{eqnarray}
est un isomorphisme. La symbole $\Thom_K(\T M)$ est une base pour le 
$\K_K(M)$-module $\K_K(\T M)$. Le caract\`ere de Riemann-Roch $RR^{K,J}(M,-)$ 
est d\'efini par le diagramme commutatif
\begin{equation}
\xymatrix@C=2cm{
 \K_K(M)\ar[r]^{\Thom_{J}} \ar[dr]_{RR^{K,J}} &
 \K_K(\T M) \ar[d]^{\indice_{M}^{K}}\\
     & R(K)\ .}
\end{equation}  

On fixe une m\'etrique riemanienne, $(-,-)_{_M}$, $K$-invariante sur $M$. 
On note $\T_K M$ le sous ensemble de $\T M$ d\'efini par
\begin{equation}
  \label{eq:def-T-K}
   \T_{K}M\ = \left\{(m,v)\in \T M,\ (v,X_{M}(m))_{_{M}}=0 \quad {\rm pour\ tout}\ 
   X\in\kgot \right\}
\end{equation}
Suivant Atiyah, un symbole $\sigma$ est dit {\em transversalement elliptique} si 
la restriction de $\sigma$ \`a $\T_{K}M$ est inversible en dehors d'un compact de  
$\T_{K}M$, i.e. $\Char(\sigma)\cap \T_{K}M$ est compact. Un symbole
transversalement elliptique $\sigma$ d\'etermine un \'el\'ement de 
$\K_{K}(\T_{K}M)$. Nous avons une application de restriction
 $\K_{K}(\T M)\to \K_{K}(\T_{K}M)$, et  Atiyah \cite{Atiyah74}  a montr\'e 
que le morphisme d'indice d\'efini sur $\K_{K}(\T M)$ s'\'etend \`a 
$\K_{K}(\T_{K}M)$ et satisfait le diagramme commutatif 
\begin{equation}\label{indice.generalise}
\xymatrix{
\K_{K}(\T M)\ar[r]\ar[d]_{\indice_{M}^K} & 
\K_{K}(\T_{K}M)\ar[d]^{\indice_{M}^K}\\
R(K)\ar[r] &  R^{-\infty}(K)\ .
   }
\end{equation}
Ici $R^{-\infty}(K)$ d\'esigne l'ensemble des caract\`eres 
g\'en\'eralis\'es de $K$: un \'el\'ement 
$\chi\in R^{-\infty}(K)$ est de la forme 
$\chi=\sum_{\mu\in\Lambda^{*}_{+}}\mm_{\mu}\, \chi_{_{\mu}}^{_K}\,$, 
o\`u $\mu\mapsto \mm_{\mu}, \Lambda^{*}_{+}\to\Z$ a une croissance au plus
polynomiale.

\medskip

{\em Th\'eor\`eme de Bott pour l'indice d'op\'erateurs homog\`enes.} Un exemple 
particuli\`erement instructif est le cas des espaces homog\`enes $M=K/H$ o\`u 
$H$ est un sous-groupe ferm\'e de $K$. 
Consid\'erons un op\'erateur diff\'erentiel $D:\Ecal^+\to \Ecal^-$ elliptique
$K$-\'equivariant sur $K/H$. Ici les fibr\'es vectoriels $K$-\'equivariants 
$\Ecal^\pm$ sont de la forme $\Ecal^\pm=K\times_H E^\pm$, o\`u $E^\pm$ sont 
 des $H$-modules. Consid\'erons la classe d\'efinie par le symbole de $D$, $\sigma(D)$,
dans $\K_{K}(\T_{K}(K/H))$. L'action de $K$ sur $K/H$ \'etant transitive,
$\T_K(K/H)$ se restreint \`a $K/H$ (vu comme la section nulle dans $\T (K/H)$).
Ainsi l'op\'erateur nul, $0_D:\Ecal^+\to \Ecal^-$ est transversalement elliptique et 
son symbole est cohomologue \`a $\sigma(D)$ dans $\K_{K}(\T_{K}(K/H))$. Si on utilise
le diagramme commutatif (\ref{indice.generalise}), on obtient alors
$$
\indice^K(D)=\indice^K(0_D)\quad {\rm dans }\quad R^{-\infty}(K).
$$
Le terme de droite de cette derni\`ere \'egalit\'e \'etant \'egal \`a
$Ker(0_D)-Coker(0_D)=L^2(\Ecal^+)-L^2(\Ecal^-)$, et comme $L^2(\Ecal^\pm)=
(L^2(G)\otimes E^\pm)^H=\indh(E^\pm)$, on obtient finalement le th\'eor\`eme de
Bott \cite{Bott65}
$$
\indice^K(D)=\indh(E^+ - E^-)\quad {\rm dans }\quad R^{-\infty}(K).
$$

\medskip

{\bf Analogies.} Pour effectuer la localisation $\K$-th\'eorique, nous travaillons de 
mani\`ere similaire que dans le contexte de la cohomologie \'equivariante: 

\begin{enumerate}
\item[$\bullet$] Le morphisme $\indice_M^K$ remplace
le morphisme d'int\'egration $\int_M$,
\item[$\bullet$] Le diagramme commutatif (\ref{indice.generalise}) remplace
$$
\xymatrix{
\Hcal_{K}^*(M)\ar[r]\ar[d]_{\int_M} & 
\Hcal_{K}^{-\infty}(M)\ar[d]^{\int_{M}}\\
S(\kgot^*)^K\ar[r] &  \fgene(\kgot)^K\ ,
}
$$
\item[$\bullet$] Le symbole $\Thom_K(M,J)$ joue le r\^ole de la classe de 
cohomologie $1_M$.
\item[$\bullet$] La $1$-forme \'equivariante $\lambda$ sera ici un champ de vecteurs 
\'equivariant.
\end{enumerate}

\medskip

Montrons comment un champ de vecteur $K$-invariant $\lambda$ r\'ealise une partition
de la base $\Thom_K(\T M)$ dans $\K_{K}(\T_{K}M)$. L'id\'ee de cette construction
est due \`a Atiyah \cite{Atiyah74}. Soit $\sigma_\lambda$ le symbole
\begin{equation}
    \sigma_\lambda(m,v):=\Thom_K(\T M )(m,v-\lambda_{m}),\quad (m,v)\in \T M.
    \label{eq:sigma.1}
\end{equation}
Le symbole $\sigma_\lambda$ est elliptique et il est homotope \`a $\Thom_K(\T M)$: ces deux
symboles d\'efinissent donc la m\^eme classe dans $\K_{K}(\T M)$. On voit que 
$\Char(\sigma_\lambda)$ est \'egal au graphe de $\lambda$, et que
$$
\Char(\sigma_\lambda)\cap\T_{K}M=\left\{(m,\lambda_{m})\in\T M, \ \ m\in 
\{\Phi_{\lambda}=0\} \right\}.
$$
o\`u $\Phi_{\lambda}:M\to\kgot^{*}$ est d\'efinie par 
$\langle\Phi_{\lambda}(m),X\rangle:= (\lambda_{m}, X_{M}\vert_{m})_{_{M}}$ 
for $X\in\kgot$.

Consid\'erons un voisinage ouvert $K$-invariant $\Ucal$ de 
$C_{\lambda}:=\{\Phi_{\lambda}=0\}$, et la restriction $\sigma_\lambda|_\Ucal$ 
qui est un symbole sur $\Ucal$. C'est un symbole transversallement elliptique car 
$\Char(\sigma_\lambda|_\Ucal)\cap\T_{K}\Ucal\simeq C_\lambda$ est compact. 
Le th\'eor\`eme d'excision donne 

\begin{lem}[\cite{pep4}]\label{lem:loc-K-theorie} On a 
$$
\Thom_K(\T M)=i_*(\sigma_\lambda|_\Ucal)\quad {\rm dans}\quad \K_{K}(\T_K M),
$$
o\`u $i:\Ucal\croc M$ est l'inclusion et $i_*:\K_{K}(\T_K \Ucal)\to\K_{K}(\T_K M)$ est 
le morphisme image directe associ\'e.
\end{lem}

Ce lemme est l'analogue du lemme \ref{lem-p-lambda} en cohomologie \'equivariante.
Dans la pratique on d\'ecompose $C_\lambda=\cup_a C^a$ en une union disjointe 
de composantes $K$-invariantes ferm\'ees, on consid\`ere des voisinages ouverts 
$K$-invariants $\Ucal^a$ de $C^a$ tels que $\Ucal^a\cap\Ucal^b=\emptyset$
si $a\neq b$. On a alors 
\begin{equation}
  \label{eq:loc-K-theorie}
  \Thom_K(\T M)=\sum_{a} i^{a}_{*}(\sigma|_{\Ucal^a})\quad 
   {\rm dans}\quad K_{K}(\T_{K}M),
\end{equation}
o\`u $i^{a}:\Ucal^a \croc M$ d\'esigne l'inclusion. Sur chaque ouvert invariant 
$\Ucal^a$, on a un morphisme canonique 
$\indice_{\Ucal^a}^{K}:K_{K}(\T_{K}\Ucal^a)\to R^{-\infty}(K)$.

\begin{defi} \label{def:RR-loc}
Pour chaque composante $C^a$ de $\{\Phi_{\lambda}=0\}$, on d\'efinit le
caract\`ere de Riemann-Roch localis\'e au voisinage de $C^a$:
\begin{eqnarray}\label{eq:RR.localise}
 RR_{C^a}^{K}(M,-) \ :\K_K(M)&\longrightarrow& R^{-\infty}(K)\\
    E\ \  &\longmapsto & 
    \indice_{\Ucal^a}^{K}(\sigma|_{\Ucal^a}\otimes p^*(E)|_{\Ucal^a} ).\nonumber   
\end{eqnarray}
\end{defi}

Le lemme \ref{lem:loc-K-theorie} donne la d\'ecomposition\footnote{Pour simplifier
la notation, on note $RR^K(M,-)$ au lieu de $RR^{K,J}(M,-)$.}
$RR^{K}(M,E)=\sum_a RR_{C^a}^{K}(M,E)$ dans $R^{-\infty}(K)$. Dans \cite{pep4}, on a 
\'etudi\'e les morphismes $RR_{C^a}^{K}(M,E)$ dans certaines situations. Les r\'esultats 
obtenus sont les analogues de ceux obtenus dans le contexte de la cohomologie
\'equivariante. Voici un bref aper\c cu des ces r\'esultats.

\bigskip

{\em Localisation de $RR^{K}(M,-)$ sur les points fixes.} 

\medskip

Dans cette partie on travaille 
avec un \'el\'ement $K$-invariant, $\beta\neq 0$, de l'alg\`ebre de Lie $\kgot$. 
Notons $\tore_\beta$ le tore de $K$ \'egal \`a 
$\overline{\{\exp(t\beta),\, t\in\R\}}$. On utilise les notations suivantes. 
Pour tout $R(K)$-module $A$, on note $A\widehat{\otimes}R(\tore_\beta)$
le $R(K)\otimes R(\tore_\beta)$-module form\'e des sommes infinies $\sum_\alpha
E_\alpha h^\alpha$ o\`u $\alpha$ parcourt l'ensemble des poids de $\tore_\beta$, et
$E_\alpha\in A$ pour tout $\alpha$.

\medskip

Consid\'erons pour l'instant le cas d'un $K$-fibr\'e vectoriel complexe 
$V\to N$, avec $N$ compact et tel que $V^\beta=N$ . Dans cette situation le symbole  
$\Thom_K(V)\in\K_K(V)$  est d\'efini de mani\`ere analogue \`a $\Thom_K(\T M)$. 
Le tir\'e en arri\`ere de $\Thom_K(V)$
sur $N$ est la classe  (d'Euler) $[\wedge^{\bullet}V]\in \K_K(N)$ qui est \'egal 
\`a la {\em diff\'erence} $[\wedge^{pair}V]-[\wedge^{impair}V]$. L'inverse 
`naturel' de $[\wedge^\bullet V]$, l'alg\`ebre sym\'etrique de $V$, ne d\'esigne pas
toujours un \'el\'ement de $\K_K(N)\widehat{\otimes}R(\tore_\beta)$. Par contre on peut 
d\'efinir un inverse $\beta$-{\em orient\'e} \cite{pep4}
\begin{equation}
  \label{eq:inverse-K-euler}
  [\wedge^{\bullet}V]^{-1}_\beta\ \in \ \K_K(N)\widehat{\otimes}R(\tore_\beta)\ .
\end{equation}
On a  $[\wedge^{\bullet}V]^{-1}_\beta= \sum_\alpha E_\alpha h^\alpha$ o\`u 
$E_\alpha\in\K_K(N)$ est non nul seulement si $\alpha=0$ ou 
$\langle\alpha,\beta\rangle> 0$. Ici la classe $[\wedge^{\bullet}V]^{-1}_\beta$
est l'analague $\K$-th\'eorique de l'inverse de la classe d'Euler \'equivariante
$\Eul{}^{-1}_\beta(V)$ dont on a parl\'e dans la premi\`ere section.
 
\bigskip

On revient au cadre d'une action de $K$ sur une vari\'et\'e $M$ munie d'une 
structure presque complexe invariante. On va maintenant d\'ecrire la localisation 
du morphisme $RR^{K}(M,-)$ que l'on obtient si l'on utilise le champ de vecteurs 
$\lambda=\beta_M$. C'est un analogue global du th\'eor\`eme de localisation
d'Atiyah-Segal-Singer pour l'indice \'equivariant. 
Comme le tore $\tore_\beta$ est dans le centre de $K$ on peut  
consid\'erer l'action de $K\times\tore_\beta$ sur $M$, et le caract\`ere
de Riemann-Roch associ\'e $RR^{K\times\tore_\beta}(M,-)$. La structure presque  
complexe sur $M$ induit une structure presque complexe sur la sous-vari\'et\'e
$M^\beta$, et une structure complexe sur le fibr\'e normal $\Ncal$ de 
$M^\beta$ dans $M$. Le caract\`ere de Riemann-Roch 
$RR^{K\times\tore_\beta}(M^\beta,-): \K_K(M^\beta)\to R(K)$ s'\'etend 
naturellement en un morphisme de $\K_K(M^\beta)\widehat{\otimes}R(\tore_\beta)$ 
dans $R(K)\widehat{\otimes}R(\tore_\beta)$.

\begin{theo}[\cite{pep4}]
Pour tout $E\in\K_K(M)$, on a
\begin{equation}
  \label{eq:loc-K-beta}
  RR^{K\times\tore_\beta}(M,E)   
=RR^{K\times\tore_\beta}\left(M^\beta,E|_{M^\beta}\otimes
[\wedge^{\bullet}\overline{\Ncal}]^{-1}_\beta\right)\quad {\rm dans}\quad 
R(K)\widehat{\otimes}R(\tore_\beta)\ .
\end{equation}
Ici $\overline{\Ncal}$ est le fibr\'e normal
muni de la structure complexe oppos\'ee.
\end{theo}

Dans mon travail concernant le calcul des $K$-multiplicit\'es de 
termes de la forme $RR^{K}(M,E)$, la formule de localisation 
(\ref{eq:loc-K-beta}) s'utilise de la mani\`ere suivante. Consid\'erons 
la multiplicit\'e de la repr\'esentation triviale $[RR^{K}(M,E)]^K$.
La th\'eor\`eme donne imm\'ediatement
$$
[RR^{K}(M,E)]^K=\Big[RR^{K\times\tore_\beta}\left(M^\beta,E|_{M^\beta}\otimes
[\wedge^{\bullet}\overline{\Ncal}]^{-1}_\beta\right)\Big]^{K\times\tore_\beta}.
$$
Dans le deuxi\`eme terme de cette \'egalit\'e on utilise maintenant le fait que 
$\tore_\beta$ agit trivialement sur $M^\beta$: pour tout fibr\'e vectoriel
\'equivariant $V\to M^\beta$, on a $[RR^{K}(M^\beta,V)]^{\tore_\beta}=
RR^{K}(M^\beta,V^\beta)$, o\`u $V^\beta$ est le sous-fibr\'e de $V$ o\`u
$\tore_\beta$ agit trivialement. Le sous-fibr\'e de 
$E|_{M^\beta}\otimes[\wedge^{\bullet}\overline{\Ncal}]^{-1}_\beta$ 
sur lequel $\tore_\beta$ agit trivialement est de rang fini: on le note 
\begin{equation}
  \label{eq:def-E-beta}
 E[\beta]= \left(E|_{M^\beta}\otimes[\wedge^{\bullet}
\overline{\Ncal}]^{-1}_\beta\right)^\beta\ .
\end{equation}
On conclut alors que 
\begin{equation}
  \label{eq:loc-beta-K-inv}
  [RR^{K}(M,E)]^K=[RR^{K}(M^\beta,E[\beta])]^K.
\end{equation}
Ainsi $[RR^{K}(M,E)]^K=0$ lorsque $E[\beta]=0$.

\bigskip
 
{\em Localisation de $RR^{K}(M,-)$ au moyen d'une application moment abstraite.}

\medskip

Dans cette section, on est toujours dans le cadre d'une action de $K$ sur 
une vari\'et\'e $M$ munie d'une structure complexe invariante. Soit 
$\Phi:M\to\kgot^*$ une {\em application moment abstraite} 
au sens de Karshon: elle est \'equivariante, et pour tout $X\in\kgot$, la 
fonction $\langle \Phi,X\rangle$ is locallement constante sur la sous-vari\'et\'e 
$M^X$. On identifie $\kgot$ \`a $\kgot^*$ au moyen d'un produit scalaire invariant.
Ici on localise le caract\`ere $RR^K(M,-)$ avec le champ de vecteurs invariant
\begin{equation}
  \label{eq:def-H}
  \lambda_m:=(\Phi(m)_M)|_m,\quad m\, \in\, M.
\end{equation}
La localisation s'effectue sur le sous-ensemble $\{\Phi_\lambda=0\}$ qui est ici 
\'egal \`a $\{\lambda=0\}$, et comme dans le cas hamiltonien on a
\begin{equation}\label{eq:beta-phi}
 \{\lambda=0\}=\bigcup_{\beta\in\Bcal_\Phi}K\left(M^\beta\cap\Phi^{-1}(\beta)\right) 
\end{equation}
o\`u $\Bcal_\Phi$ est un sous-ensemble fini d'une chambre de Weyl $\tgot^*_+$.
D'apr\`es la d\'efinition \ref{def:RR-loc}, on a la d\'ecomposition
\begin{equation}
  \label{eq:decomposition-K}
  RR^K(M,-)=\sum_{\beta\in\Bcal_\Phi}RR^K_\beta(M,-)\ ,
\end{equation}
o\`u $RR^K_\beta(M,-):\K_K(M)\to R^{-\infty}(K)$ est le caract\`ere de Riemann-Roch  
localis\'e au voisinage de $K\left(M^\beta\cap\Phi^{-1}(\beta)\right)$ (au moyen 
du champ de vecteur $\lambda$). Pour tout $E\in\K_K(M)$, la multiplicit\'e 
$[RR^{K}(M,E)]^K$ v\'erifie $[RR^{K}(M,E)]^K=\sum_\beta [RR^{K}_\beta(M,E)]^K$. 
Dans la pratique je calcule $RR^{K}_0(M,E)$ et je montre que dans certains 
cas \break $[RR^{K}_\beta(M,E)]^K=0$ pour tout $\beta\neq 0$. Je termine
cette section en donnant un bref aper\c cu des r\'esultats obtenus sur ces 
$RR^K_\beta(M,-)$.

\bigskip

\underline{Le cas $\beta=0$.} 

\bigskip

Lorsque $0$ est une valeur r\'eguli\`ere de $\Phi$, je montre dans \cite{pep4} que 
la structure presque complexe $K$-invariante sur $M$ induit une 
structure $\spinc$ sur la $V$-vari\'et\'e $\Mcal_0:=\Phi^{-1}(0)/K$. Soit  
$\Qcal(\Mcal_0,-): \K(\Mcal_0)\to \Z$ le morphisme d\'etermin\'e par cette 
structure $\spinc$. Dans ce cas, je calcule $RR^K_0(M,-)$. Je montre en particulier 
que la multiplicit\'e de la repr\'esentation triviale dans $RR^K_0(M,E)$ est 
$\Qcal(\Mcal_0,\Ecal_0)$ avec $\Ecal_0:=(E|_{\Phi^{-1}(0)})/K$. Un calcul similaire
avait \'et\'e \'effectu\'e par Vergne dans le cas d'une action hamiltonienne 
du cercle \cite{Vergne96}.

\bigskip

\underline{Le cas $\beta\neq 0$ central.}

\bigskip

Ici on consid\`ere la $K$-vari\'et\'e des points fixes $M^\beta$, munie
de la structure presque complexe induite et de l'application moment abstraite 
$\Phi|_{M^\beta}$. On peut alors d\'efinir le caract\`ere de Riemann-Roch 
$RR^K_\beta(M^\beta,-)$ localis\'e au voisinage de $M^\beta\cap\Phi^{-1}(\beta)$. 
Comme pour la localisation sur les points fixes on \'etend ce morphisme en
\break $RR^{K\times\tore_\beta}_\beta(M^\beta,-)$. Je montre alors  que
\begin{equation}
  \label{eq:loc-beta-gene}
   RR^{K\times\tore_\beta}_\beta(M,E)   
=RR^{K\times\tore_\beta}_\beta\left(M^\beta,E|_{M^\beta}\otimes
[\wedge^{\bullet}\overline{\Ncal}]^{-1}_\beta\right)
\end{equation}
dans  $R^{-\infty}(K)\widehat{\otimes}R(\tore_\beta)$. Dans cette situation 
on a un raffinement de (\ref{eq:loc-beta-K-inv}) 
$$
[RR^{K}_\beta(M,E)]^K=[RR^{K}_\beta(M^\beta,E[\beta])]^K,
$$
o\`u $E[\beta]$ est le sous-fibr\'e de
$E|_{M^\beta}\otimes[\wedge^{\bullet}\overline{\Ncal}]^{-1}_\beta$ 
sur lequel $\tore_\beta$ agit trivialement. Lorsque $E$ est un fibr\'e strictement
$\Phi$-positif (voir d\'efinition \ref{def:stric-positif}) on a $E[\beta]=0$, et 
donc $[RR^{K}_\beta(M,E)]^K=0$.

\bigskip

\underline{Le cas $K_\beta\neq K$.}

\bigskip

Dans cette situation, j'obtiens des formules d'induction qui me permettent
de `revenir' au cas pr\'ec\'edent. Le choix d'une chambre de Weyl d\'etermine 
naturellement une structure structure complexe sur $\kgot/\kgot_\beta$. 

Nous avons le morphisme d'induction $\indB:\fgene(K_\beta)^{K_\beta}\to
\fgene(K)^{K}$ o\`u $\fgene(K_\beta)$ est l'ensemble des fonctions g\'en\'eralis\'ees
sur $K_\beta$, et les invariants sont pris par rapport \`a l'action de conjugaison.
L'application $\indB$ est d\'etermin\'e par la relation: pour tout 
$\phi\in\fgene(K_\beta)^{K_\beta}$, on a $\int_K\indB(\phi)(k)f(k)dk=cst
\int_{K_\beta}\phi(h)f(h)dh$ pour tout $f\in\f(K)^K$,  o\`u $cst=\Vol(K,dg)/
\Vol(K_\beta,dh)$. Nous avons aussi le morphisme d'induction 
holomorphe $\HolB:R^{-\infty}(K_\beta)\to R^{-\infty}(K)$ d\'etermin\'e par la relation:
$\HolB(\chi)=\indB\left(\chi\wedge_\C^\bullet\kgot/\kgot_\beta\right)$ pour tout 
$\chi\in R^{-\infty}(K_\beta)$.

Nous pouvons consid\'erer la vari\'et\'e $M$ munie de l'action du sous-groupe 
$K_\beta$, et de l'application moment abstraite $\Phi_{K_\beta}:M\to\kgot_\beta^*$
qui est le compos\'e de $\Phi$ avec la projection $\kgot^*\to\kgot^*_\beta$.
Nous d\'efinissons alors le caract\`ere de Riemann-Roch 
$RR^{K_{\beta}}_{\beta}(M,-)$ localis\'e au voisinage de 
$M^\beta\cap(\Phi_{K_\beta})^{-1}(\beta)$. Je montre alors que pour tout
$E\in\K_K(M)$ on a

\begin{equation}\label{eq.relation.induction.1}
 RR^{K}_{\beta}(M,E)=\HolB \left(RR^{K_{\beta}}_{\beta}(M,E)\,
 \wedge_{\C}^{\bullet}\overline{\kgot/\kgot_{\beta}}\right).
\end{equation}

Cette formule d'induction combin\'ee avec (\ref{eq:loc-beta-gene})
permet de montrer que \break $[RR^{K}_{\beta}(M,E^{\otimes l})]^K=0$ pour $l$ assez grand, 
lorsque $E$ est un fibr\'e strictement $\Phi$-positif.

\bigskip 

Lorsque l'application $\Phi$ est une (vraie) application moment associ\'ee \`a
une {\em action hamiltonienne} de $K$ sur $(M,\omega)$, je raffine 
(\ref{eq.relation.induction.1}) de la fa\c con suivante. Consid\'erons la face 
ouverte $s$ de la chambre de Weyl qui contient $\beta$ et la `tranche
symplectique' $\Ycal_s$ associ\'ee. C'est une sous-vari\'et\'e symplectique de 
$M$ munie d'une action hamiltonienne de $K_\beta$: notons 
$\Phi_s:\Ycal_s\to\kgot_\beta^*$ l'application 
moment correspondante. On munit les vari\'et\'es symplectiques $M$ et
$\Ycal_s$ de structures presque complexe compatibles. Sur $\Ycal_s$ on consid\`ere
le caract\`ere de Riemann-Roch 
$RR^{K_{\beta}}_{\beta}(\Ycal_s,-)$ localis\'e au voisinage de 
$\Ycal_s^\beta\cap(\Phi_s)^{-1}(\beta)=M^\beta\cap\Phi^{-1}(\beta)$. Je montre 
 alors que pour tout $E\in\K_K(M)$ on a

\begin{equation}\label{eq.relation.induction.2}
RR^{K}_{\beta}(M,E)=
\HolB\left(RR^{K_{\beta}}_{\beta}(\Ycal_{s},
E\vert_{\Ycal_{s}})\right)\ .
\end{equation}

L'expression (\ref{eq.relation.induction.2}) permet de montrer que 
$[RR^{K}_{\beta}(M,E)]^K=0$ dans l'un des deux cas suivants:
\begin{enumerate}
\item[$\bullet$] $E$ est strictement $\Phi$-positif et $0\notin\Phi(M)$,
\item[$\bullet$] $E$ est $\Phi$-positif et $0\in\Phi(M)$.
\end{enumerate}


\end{document}